\documentclass[final]{siamltex}

\usepackage{latexsym, enumerate}
\usepackage{eepic}
\usepackage{epic}
\usepackage{graphicx}
\usepackage{color}
\usepackage{ifpdf}
\usepackage{amssymb,amsmath,epsfig, algorithm,algorithmicx,multirow}
\usepackage{amsfonts, dsfont}
\usepackage{subfigure}
\usepackage{multirow}
\usepackage{makecell}
\usepackage{bm}
\usepackage{mathrsfs}
\usepackage{relsize}

\newcommand{\bff}{\mathbf}
\newcommand{\bb}{\mathbb}
\newcommand{\ca}{\mathcal}

\newtheorem{thm}{Theorem}[section]

\newtheorem{rem}{Remark}[section]

\title{Bayesian identification of discontinuous fields with an ensemble-based variable separation multiscale method
\thanks{L. Jiang acknowledges the support of Chinese NSF 11471107.
}}

\author{Na Ou\thanks{College of Mathematics and Econometrics, Hunan University, Changsha 410082, China.({\tt oyoungla@hnu.edu.cn}).}
\and
Guang Lin\thanks{Department of Mathematics $\&$ School of Mechanical Engineering, Purdue University, West Lafayete, IN 47907-2067, USA.({\tt guanglin@purdue.edu}).}
\and
Lijian Jiang\thanks{School of mathematical sciences, Tongji University, Shanghai, 200092, China.({\tt ljjiang@tongji.edu.cn}), Corresponding author.}
}
\begin{document}

\maketitle

\begin{abstract}
This work presents a multiscale model reduction approach to discontinuous fields identification problems in the framework of Bayesian inference. An  ensemble-based variable separation (VS) method is proposed to approximate multiscale basis functions used to build  a coarse model.  The variable-separation  expression is constructed  for stochastic multiscale basis functions based on the random  field, which is treated Gauss process as prior information. To this end,  multiple local inhomogeneous Dirichlet boundary condition problems are required to be solved, and the ensemble-based method is used to obtain variable separation forms for the corresponding local functions. The local functions share the same interpolate rule for different physical basis functions  in each coarse block. This approach significantly  improves the efficiency of computation. We obtain the variable separation expression of multiscale basis functions, which  can be used  to the models with different boundary conditions and source terms, once the expression constructed. The proposed method is applied to discontinuous field identification problems where the hybrid of total variation and Gaussian (TG) densities are imposed as the penalty. We give  a  convergence analysis of the approximate posterior to the reference one with respect to the Kullback-Leibler (KL) divergence under the hybrid prior. The  proposed method is applied  to identify  discontinuous structures in permeability  fields.
Two patterns of discontinuous structures are considered in numerical examples: separated blocks and nested blocks.
\end{abstract}



\pagestyle{myheadings}

\thispagestyle{plain}
\markboth{Na Ou, Guang Lin and L. Jiang}{Bayesian discontinuous fields identification}

\section{Introduction}

The characterization of the permeability is critical  for the modeling subsurface  flow and transportation. Properly quantifying the uncertainties induced by the permeability field is
very important for making reliable probabilistic based predictions and future decisions. The involved unknown spatially distributed parameters are need to be estimated  in some quantified sense, such that the observations can be best explained.

Because  the dimension of measurement data is usually much smaller than the dimension of the unknown spatially distributed parameters, the limited number of measurements can only depict the permeability field with high uncertainty. The weak sensitivity of measurements to the unknowns and the noise in the observations, often results in ill-posedness of the inverse problems. The prior regularizing information is naturally utilized as a penalty to remedy this issue. Gaussian measures are widely used as the prior distributions, which also work as the Tikhonov regularization term used from the view of optimization problems. A rigorous Bayesian framework for the inverse problems in function spaces  is developed in \cite{stuart2010inverse} for smooth unknown random fields.

However, in practical situations, the spatially varying fields can be better represented as consisting of a few relatively uniform geologic facies or zones with abrupt changes at their boundaries.
One approach to facies detection is the level-set method \cite{osher1988fronts, burger2001level, dunlop2017hierarchical}, in which the shapes of facies are not required to be predetermined, to identify zonal structures. The geometry of the discontinuity of the permeability is represented implicitly by level set functions. The authors of \cite{chan2004level, chung2005electrical} combine the level-set method with the total variational (TV) norm \cite{rudin1992nonlinear} to control both the jumps in the permeability and the length of the level sets. The TV norm is also used as the prior density in \cite{lee2013bayesian} to describe the sharp changes of the spatially varying field in the framework of Bayesian inference. The well-defined posterior measure in the function space is developed in paper \cite{yao2016tv}, where TV-Gaussian (TG) prior is presented to address the infinite-dimensional Bayesian inverse problems, for discontinuous random fields. We choose the TG prior to identify the discontinous permeability field in this paper.

Markov chain Monte Carlo (MCMC) methods \cite{liu2008monte, brooks2011handbook, bui2013computational} are the standard techniques for sampling the posterior distribution in Bayesian inferences. An infinite number of samples are required to pursue the convergence and reliability of the chain, which implies we need to simulate large-scale PDE-based forward models hundreds of thousands or even millions times, for characterizing the posterior distribution. The discritization of the spatial domain results in high-dimensional parameters inevitably, while the efficiency of many standard MCMC chains degrades as the dimension of the parameters increases. The computational defect of standard MCMC method has promoted the development of efficient MCMC techniques, such as incorporating the data information to improve the sampling efficiency: the two-stage and cheaper approximation strategies \cite{efendiev2005efficient, christen2005markov}, the stochastic Newton MCMC \cite{martin2012stochastic} and the dimension-independent likelihood-informed MCMC \cite{cui2016dimension}. The preconditioned Crank-Nicolson MCMC (pCN-MCMC) method described in \cite{cotter2013mcmc}, of which the acceptance probability is invariant with respect to the dimension of the parameter space, is applied to draw samples from the Bayesian posterior in the present work.

As noted before, exploration of the posterior distribution requires repeated evaluations of the forward operator, and a  large number of samples are needed to ensure reliable estimates of the inference. One attempts at accelerating Bayesian inference in computationally intensive inverse problems have relied on reduction of the stochastic forward model \cite{frangos2010surrogate}, the approaches of building stochastic surrogates include generalized polynomial chaos (gPC)-based stochastic method \cite{xiu2002wiener, marzouk2009stochastic, zhou2015weighted, yan2015stochastic, jiang2018bayesian}, Gaussian process \cite{ras2006gaussian, bilionis2013solution} or projection-type reduced order models \cite{hinze2005proper, noor1980reduced, nouy2010proper, li2017novel}, for multiple solutions, proper orthogonal decomposition (POD) \cite{hinze2005proper} reduced-order model is incorporated into the ensemble-based method \cite{jiang2014algorithm} to significantly reduce the costs \cite{gunzburger2017ensemble}, etc. Large numbers of forward model simulations are required for these methods, especially when the dimension of stochastic parameter is high. The multi-fidelity model \cite{Robinson2008Surrogate, Zhu2014Computational} technique is one of the method in reducing the computational time, it is widely used in uncertainty propagation, inference, and optimization \cite{forrester2007multi, peherstorfer2018survey} to improve the efficiency, e.g., the author of \cite{yan2019adaptive} use the adaptive multi-fidelity strategy to obtain the corrected gPC approximation, and apply the surrogate in nonlinear Bayesian inverse problems. Multiscale Finite Element Method (MsFEM) \cite{hou1997multiscale} is one of the model reduction methods and gives a coarse-scale equation  through variational formulation  by multiscale basis functions, which contain fine scale information. The generalized Multiscale Finite Element Method (GMsFEM) \cite{efendiev2011multiscale, efendiev2013generalized} inherits  the important features from MsFEM, once the multiscale basis functions constructed, it can be used repeatedly for models with different source terms and boundary conditions, while the coarse space in GMsFEM is more flexible.

In the paper, we develop an  ensemble-based variable separation multiscale (VS-Ms) method in accelerating the forward model simulation. First, we build local problems based on  GMsFEM, with multiple determined inhomogeneous Dirichlet boundary conditions, to obtain the reduced snapshot space. We use the VS method to solve the local problems in each coarse neighborhood such that the local functions can be expressed in variable separation form. However, it would be time consuming if we approximate each local function by the VS procedure. On the other hand, as the number of multiscale basis functions retained in the coarse block increase, the time for calling the variable separation form may increase to a large magnitude, even the separated terms is fixed at a small number.

To efficiently construct multiscale basis functions in random  space, the ensemble-based method is used in this paper. We run VS process to solve the local problems with the ensemble of the inhomogeneous Dirichlet boundary conditions, and then reverse the physical basis for each local function basing on the mechanism and residual equation of the VS process. Then we have stochastic multiscale basis functions in one coarse block sharing  the same interpolate rule, i.e., for any random parameter, the separated stochastic functions from the VS form need only to be called one time, and the multiscale basis functions can then be approximated by multiplying the calculated stochastic part with the corresponding physical parts. The derived  multiscale basis functions can capture most characteristics of the ones built by directly solving local problems.

The proposed method shares  the merit from the GMsFEM and can be applied to the model with different boundary conditions and source terms, once the stochastic multiscale basis functions expression is constructed. We apply it to discontinuous field identification  problems under the framework of Bayesian inference. The variable separated form of the multiscale basis functions are constructed with respect to Gaussian random field, of which the covariance matrix depends on the relative distance of the spatial points. In practice, the ensemble-based VS-MsFEM approach is used only $4$ times in constructing variable separated form of multiscale basis functions. More details on the proposed method will be displayed in Section $\ref{En-VS-Ms}$.
The number of local basis functions selected in each coarse neighborhood dominates the effects on the accuracy.  The proposed method is applied in recovering discontinuous fields with different structures, e.g., two separated blocks and nested blocks. We give a brief convergence analysis of the approximate posterior to the reference one with respect to the Kullback-Leibler divergence, and the numerical results  verify the convergence  of the proposed method.

The outline of this paper is organized as follows. In Section $2$,  we briefly introduce   the Bayesian inference with hybrid prior, the pCN-MCMC algorithm and the general convergence result  of the approximated posterior, which is induced by the reduced model, under the TG prior. Section $3$ describes the proposed model reduction method, where the local basis functions stem from the snapshots of  GMsFEM, and the VS method is applied to inhomogeneous Dirichlet boundary conditions and the ensemble-based VS method is used to obtain the variable expression for the problems with multiple boundary conditions. Some numerical examples are presented in Section $4$, we carry out a
few numerical  simulations to verify the efficiency and accuracy of our proposed method. The proposed approach is also applied to posterior densities approximation, where discontinuous fields with different structures are identified. Some conclusions and comments are made finally.

\section{Problem setting}
We consider a system described by partial differential equations
\begin{equation}\label{PDE}
\ca{F}(u; \kappa(x))=f, \quad x \in D \subset \bb{R}^2,
\end{equation}
where the unknown discontinuous coefficient $\kappa(x) \in L^{\infty,+}(D)$,
\[
 L^{\infty,+}(D):=\{L^\infty(D)|\kappa(x) > 0\}.
\]
The operator $\ca{F}$ represents the forward model, which describes the relation between the coefficient  $\kappa(x)$, external source/sink term $f$ and the output state $u$. The coefficient identification problem for $(\ref{PDE})$ seeks to determine the coefficient $\kappa(x)$ in such a way that the observed output matches the indirect measurement $d$ in a prescribed sense, the data $d$ is related to $\kappa(x)$ via
\begin{equation}\label{noise}
d={\bff G}(\kappa)+\zeta,
\end{equation}
where $\bff G: L^{\infty,+}(D)\rightarrow \bb{R}^{n_d}$ and $\zeta$ is an $n_d$-dimensional independent and identically distributed (i.i.d.) Gaussian random vector with mean zero and standard deviation $\sigma$. The likelihood is given by
\begin{equation}
\label{LH}
\pi(d|\kappa) = (2\pi\sigma^2)^{-\frac{n_d}{2}}\exp\bigg(-\frac{\|d-{\bff G}(\kappa)\|_2^2}{2\sigma^2}\bigg).
\end{equation}
Due to the sparsity of the given measurement, some regularization terms are required. In order to guarantee the well-posedness of the inverse problem, the weighted $L_2$ norm penalty is imposed.  We set
\begin{equation}
\label{GP}
\log{\kappa(x)} \sim {\mathcal {GP}}(0, \ca{C}(x_1, x_2)) ,
\end{equation}
where $\ca{C}(\cdot, \cdot)$ is the kernel operator, ${\mathcal {GP}}(\cdot, \cdot)$ represents the Gaussian process. We concentrate on the  squared exponential kernel in this paper, i.e., it has the form
\[
\ca{C}(x_1, x_2)=\sigma_x^2\exp\bigg(-\frac{(x^1_1-x^1_2)^2}{2l_{x_1}^2}-\frac{(x^2_1-x^2_2)^2}{2l_{x_2}^2}\bigg),
\]
where $\sigma_x^2$ is the variance of the random field, $l_{x_1}$ and $l_{x_2}$ are the length scales along the direction $x_1$ and $x_2$, respectively. Other penalty term is also required to pursue the discontinuous property of coefficient. We follow \cite{yao2016tv} and the edge pursuing penalty term $\mathcal{R}(\kappa)$ is set as total variation norm, i.e.,
\[
\mathcal{R}(\kappa)=\lambda\|\kappa(x)\|_\text{TV}=\lambda \int \|\nabla\kappa(x)\|_2 dx,
\]
where $\lambda$ is a prescribed positive constant.
The infinite problem has been proved well-defined in paper \cite{yao2016tv}, we discretize the spatially varying field matching with the grid, on which the forward model is solved, i.e., $\kappa(x)$ can be approximated by
\begin{equation}\label{affinem}
\kappa(x):=\kappa(x; \xi)=\sum_i^{N_{x_1}N_{x_2}} \mathbf{I}_{\{x\in D_i\}}{\text e}^{\xi_i},
\end{equation}
where $N_{x_1}$ and $N_{x_2}$ are the number of partition nodes in the $x_1$ and $x_2$ direction, respectively. $\bigcup_{i} D_i$ are the non-overlapping partition of the domain $D$, $\mathbf{I}_{\{x\in D_i\}}$ is the indicator function, and
\begin{equation}\label{dGP}
\xi\sim {\mathcal N}(0, \Sigma),
\end{equation}
where $\Sigma$ is the covariance of the joint Gaussian distribution. Then we make inference to the unknown $\xi \in {\mathbb R}^{n_\xi}$, where $n_\xi=N_{x_1}N_{x_2}$ is the dimension of $\xi$.  Then the posterior expression is
\begin{equation}
\label{Bayesxi}
\pi(\xi|d)\propto \exp\bigg(-\frac{\|d-{\bff G}(\xi)\|_2^2}{2\sigma^2}-\frac{1}{2}\|\xi\|_{\Sigma}^2-\mathbf{R}(\xi)\bigg),
\end{equation}
where
\[
\|\xi\|_{\Sigma}^2=\|\Sigma^{-\frac{1}{2}}\xi\|_2^2,
\]
and $\mathbf{R}(\xi)$ is the discretized approximation with respect to $\xi$,
more details about the discretization of the TV norm can be referred to \cite{vogel2002computational}(chapter 8.2). Posterior samplers are obtained by pCN-MCMC scheme \cite{cotter2013mcmc}.



We use the ensemble-based VS multiscale (VS-MS) method, which will be described in detail in Section 3, to accelerate the sampling. Let $\bff {G_N}$ be the approximated forward observation operator. Then the corresponding approximated posterior is
\begin{equation}\label{aBayesxi}
\pi_{\bff{N}}(\xi|d)\propto \exp\bigg(-\frac{\|d-\bff {G_N}(\xi)\|_2^2}{2\sigma^2}-\frac{1}{2}\|\xi\|_{\Sigma}^2-\mathbf{R}(\xi)\bigg).
\end{equation}
For the convenience of notation, we will use $\pi^d(\xi)$ and $\pi_{\bff{N}}^d(\xi)$ to denote the posterior density $\pi(\xi|d)$ and $\pi_{\bff{N}}(\xi|d)$, respectively.

We use Kullback-Leibler (KL) divergence \cite{marzouk2009dimensionality} to measure  the difference between the approximated posterior and the reference one. For probability density functions $\pi_{\bff{N}}^d(\xi)$ and $\pi^d(\xi)$, KL divergence is defined by
\[
D_{KL}(\pi_{\bff{N}}^d||\pi^d)=\int \pi_{\bff{N}}^d(\xi) \log \frac{\pi_{\bff{N}}^d(\xi)}{\pi^d(\xi)}d\xi.
\]
$D_{KL}$ measures the difference between the two probability distributions and is nonnegative. For a given discretization, we have a brief convergence analysis without regard to the discrete error of the unknowns and the total variation penalty.
\begin{thm}\label{DKL}
(Theorem 3.1 of \cite{yan2017convergence}.)Suppose the functions $\bff G$ and $\bff {G_N}$ are under some assumption, and the observational error has an i.i.d. Gaussian distribution. Then the approximation posterior $\pi_{\bff N}^d$ and the true posterior density $\pi^d$ are close with respect to the Kullback-Leibler distance, there is a constant $C$ independent of $\bff N$, such that
  \[
  D_{KL}(\pi_{\bff N}^{d}\|\pi^d)\leq \frac{C}{\sigma^4}\|\bff G(\xi)-\bff {G_N}(\xi)\|_{L_{\pi_0}^2}^2.
  \]
\end{thm}

\begin{rem}
The assumptions functions $\bff G$ and $\bff {G_N}$ satisfied in Theorem \ref{DKL} are: the forward operator $\bff G$ satisfies $\sup_\xi \| \bff G (\xi)\|_2<\infty$, which can be met in many applications. The surrogate $\bff {G_N}$ satisfies $\sup_\xi \| \bff G (\xi)-\bff {G_N}(\xi)\|_2\rightarrow 0$ as $\bff N\rightarrow \infty$. Here $\bff N$ mainly represents the number of multiscale basis functions on each coarse block. We note the error term $\|\bff G(\xi)-\bff {G_N}(\xi)\|_{L_{\pi_0}^2}^2$ is calculated under the Gaussian prior, but not the hybrid TG prior, this is due to the non-negative and boundedness of the term ${\mathbf R}(\xi)$, it can be treated similarly to the forward model, hence the consequence of Theorem 3.1 in paper \cite{yan2017convergence} is still applicable.
\end{rem}

%
%
%
%

\section{Ensemble-based VS-MS method}
\label{En-VS-Ms}
The dimension of $\xi$ depends on the discretization of the physical domain $D$, which could be very high. If constructing the surrogate model directly, it will suffers from severe curse of dimensionality. For a change, we construct surrogate model over local domains. We construct local problems stemming from  GMsFEM, and then approximate the  local functions on each coarse neighborhood by the ensemble-based VS scheme. The resulted variable separation expressions for stochastic multiscale basis functions share the common interpolate rule, which save the computational cost further.

   We use the following model equation subject a boundary condition to show the construction of ensemble-based VS-MS basis functions,
   \[
   -\text{div} (\kappa(x; \xi)\nabla u)=f(x) \quad \text{in} \quad D.
   \]
   This equation can model the single-phase subsurface flow. We use generalized   finite element method  \cite{ bbo03}  to solve this equation on a coarse grid, where the basis functions are constructed by the proposed  ensemble-based VS method.
   We refer this method to an ensemble-based VS-MS method.  The generalized   finite element method is standard (see \cite{ bbo03}).  In the paper, we will focus on the construction of multiscale basis functions on each coarse block.

\subsection{Local problems based on GMsFEM}

\begin{figure}
  \centering
  \includegraphics[width=0.6\textwidth]{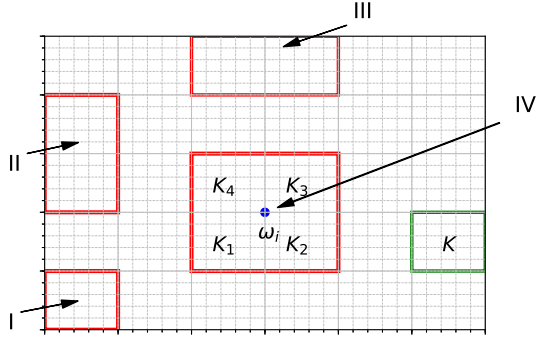}
  \caption{Illustration of the 4 types of coarse neighborhood and the coarse element.}
   \label{coarse-cell}
\end{figure}

For GMsFEM, we need to pre-compute a set of multiscale basis functions. For a random sample taken from the prior $(\ref{GP})$, we solve the following local problem on each coarse block $\omega_i$ to construct the snapshot space,
\begin{eqnarray}
\label{snapshots}
\begin{cases}
& -\text{div}(\kappa(x; \xi)\nabla \varphi^i_l) =0 \ \ \text{in}\ \omega_i\\
&\varphi^i_l = \delta^h_l(x)\ \ \text{on}\ \partial \omega_i,
\end{cases}
\end{eqnarray}
where functions $\delta^h_l(x)$ are defined with respect to the fine grid boundary point on $\partial \omega_i$. The space of snapshots $V^{\omega_i}_{\text{snap}}$ can be constructed  as
\[
V^{\omega_i}_{\text{snap}}=\text{span}\{\varphi^i_l, 1 \leq i \leq N_H, 1 \leq l \leq M^i_{\text{snap}}\},
\]
where $N_H$ is the number of coarse nodes and $M^i_{\text{snap}}$ is the number of boundary nodes on $\partial \omega_i$. The snapshot functions can be stacked into a matrix as
\[
R_\text{snap}=[\varphi_1,\cdots,\varphi_\text{Msnap}].
\]
In order to reduce the dimension of the snapshot space, we perform a spectral decomposition of the space of snapshots and solve the local problems
  \begin{eqnarray}
  \label{redsnapshots}
   \begin{cases}
   & -\text{div}(\kappa(x; \xi)\nabla \psi^i_j) =\lambda_j \kappa(x; \xi) \psi^i_j\ \ \text{in}\ \omega_i\\
   & \kappa(x; \xi)\nabla \psi^i_j\cdot \vec{n} = 0\ \ \text{on}\ \partial \omega_i,
  \end{cases}
  \end{eqnarray}
where $\psi_j^i \in V^{\omega_i}_{\text{snap}}$. We choose the smallest $M_i$ eigenvalues of
  \[
  \tilde A\psi^i_j=\lambda_j \tilde S\psi^i_j,
  \]
and take the corresponding eigenvectors in the space of snapshots by setting $\psi_j^i=\sum_l \psi_{j,l}^i\varphi_l$, for $j=1,\cdots,M_i$, to form the reduced snapshot space, where $\psi_{j,l}^i$ are the coordinates of the vector $\psi^i_j$ and
  \[
  \tilde A=[a_{mn}]=\int_{\omega_i} \kappa(x; \xi)\nabla \varphi_n\nabla \varphi_m=R_{\text{snap}}^TAR_{\text{snap}},
  \]
  \[
  \tilde S=[s_{mn}]=\int_{\omega_i} \kappa(x; \xi) \varphi_n\varphi_m=R_{\text{snap}}^TSR_{\text{snap}}.
  \]
where $A$ and $S$ denote fine-scale matrices corresponding to the stiffness and mass matrices with the coefficient  $\kappa(x; \xi)$.  The relationship between a coarse neighborhood and its coarse elements is illustrated in Figure \ref{coarse-cell}.

The goal of Eqs.$(\ref{snapshots})$ and $(\ref{redsnapshots})$ is to construct a few effective boundary conditions in each coarse block, so that the functions $\psi^i_j(x, \xi), j=1, \cdots, M_i$ expanding the reduced snapshot space can be used to approximate the solution space. The construction of the reduced snapshot space is equivalent to solving the equation directly
\begin{eqnarray*}
\begin{cases}
& -\text{div}(\kappa(x; \xi)\nabla \phi^i_j) =0, \ \ \text{in}\ \omega_i,\\
&\phi^i_j = \psi^i_j(x, \xi), \ \ \text{on}\ \partial \omega_i
\end{cases}
\end{eqnarray*}
in the fine grid finite element space, where $\psi^i_j(x; \xi)$ is treated as the effective boundary condition. However, it is inadvisable to build the reduced snapshot space by solving the equation above, since $\psi^i_j(x, \xi)$ depends on the stochastic variable  $\xi$, and we have no explicit expression of it. As a compromise, we solve local problems
\begin{eqnarray}
\label{effMS-basis}
\begin{cases}
& -\text{div}(\kappa(x; \xi)\nabla \phi^i_j) =0, \ \ \text{in}\ \omega_i,\\
&\phi^i_j = \bar{\phi}^i_j(x), \ \ \text{on}\ \partial \omega_i
\end{cases}
\end{eqnarray}
to construct the reduced snapshot space, where $\bar{\phi}^i_j(x)=\psi^i_j(x, \bar{\xi})$ is built with $\kappa(x; \xi)=\kappa(x; \bar{\xi})$. Here $\bar{\xi}$ is an arbitrary realization taken from the prior, e.g., the mean of the random field. In this way, the boundary condition $\bar{\phi}^i_j(x)$ is deterministic and pre-computed  for any  $\xi$, then the local problem $(\ref{effMS-basis})$ can be solved directly.

It is reasonable to treat $\bar{\phi}^i_j(x)$ as the effective boundary condition, since it is constructed based on the model information. Even when the model information is rough, the boundary conditions can sometimes be constructed effectively, which further leads to good approximation of the solution.
Finally, we multiply the partition of unity functions (e.g., linear hat functions over coarse elements) by $\phi^i_j(x, \xi)$ obtained through solving $(\ref{effMS-basis})$ to construct the local basis space.

\subsection{VS method}
\label{VS}
We use the VS method proposed in \cite{li2017novel} to construct surrogate model for solutions to Eq.$(\ref{effMS-basis})$. For the convenience of notification, we leave out the scripts in $\phi^i_j(x, \xi)$, $\bar{\phi}^i_j(x)$ and $\omega_i$
in this section, to demonstrate the application of VS method to inhomogeneous Dirichlet boundary condition problems.
Let $H^1_0$ be the usual Sobolev space.  A weak formulation of Eq. $(\ref{effMS-basis})$ can be written as: find $\phi^0 \in \ca{V}$ such that
\[
a(\phi^0(\xi), v; \xi)=-a(\bar{\phi}, v; \xi),  \ \ \forall v\in \ca{V} ,
\]
where $\ca{V}:=H^1_0(\omega)$, $a(\cdot ; \cdot)$ is a bilinear form on $\ca{V}$ from the diffusion equation, and $\phi=\phi^0+\bar{\phi}$ is the solution. The bilinear form $a(\cdot, \cdot; \xi)$ is assumed to be affine with respect to $\xi$ i.e.,
\begin{equation}\label{affine}
a(w, v; \xi) =\sum_{p=1}^{m_a}\kappa_p(\xi)a_p(w, v), \ \forall v,w \in {\ca V}, \ \forall \xi \sim {\mathcal N}(0, \Sigma),\\
\end{equation}
{where $\kappa_p(\xi)$ is stochastic functions with respect to $\xi$, $a_p : \ca{V} \times \ca{V} \longrightarrow \mathbb{R}$ is a symmetric bilinear form. If $a(\cdot, \cdot; \xi)$ is not affine, the VS method \cite{li2017novel} can still be used. Let $\ca{V}_h \subset \ca{V}$ be given finite dimensional approximation space. We seek the solution of Eq.$(\ref{effMS-basis})$ approximated by
\begin{equation*}
\phi(x;\xi)\approx \phi_N^0(x;\xi)+\bar{\phi}(x),
\end{equation*}
and $\phi_N^0(x;\xi)$ has the separated form
\begin{equation}
\label{vs-form}
\phi^0_N(x;\xi)=\sum_{q=1}^N \Phi_q(x)\eta_q(\xi),
\end{equation}
where $\eta_q(\xi)$ are stochastic functions and $\Phi_q(x) \in \ca{V}_h$ are deterministic functions, the separated term $N$ may be much less than the degree of freedom of the original problem.

Let the residual for the VS procedure be
\[
e(\xi):=\phi^0(\xi)-\phi^0_{k-1}(\xi).
\]
Then,  for $\forall  v\in \ca{V}_h$,
\[
a(e(\xi), v; \xi)=-a(\bar{\phi}, v; \xi)-a(\phi^0_{k-1}(\xi), v; \xi).
\]
Let $r(v; \xi)\in \ca{V}_h^*$ be the residual
\begin{equation}\label{residual}
r(v; \xi)=\begin{cases}
-a(\bar{\phi}, v; \xi), & k=1, \\
-a(\bar{\phi}, v; \xi)-a(\phi^0_{k-1}(\xi), v; \xi), & k\geq 2.
\end{cases}
\end{equation}
We have
\begin{equation}
\label{e-solution}
a(e(\xi), v; \xi)=r(v; \xi), \ \forall  v\in \ca{V}_h.
\end{equation}

Let $\Xi$ be  a small sample set taken from ${\mathcal N}(0, \Sigma)$. At step $k$, we choose
\[
\xi^k:=\begin{cases}
\text{randomly chosen from}\ \ \Xi, \ &k=1,\\
\arg\max_{\xi\in\Xi}\|\hat e(\xi)\|^2_{\ca V}, \ & k\geq2.
 \end{cases}
 \]
where the error indicator $\|\hat e(\xi)\|^2_{\ca V}$ has the form
\begin{equation}
\label{e-hat-norm}
\begin{split}
\|\hat e(\xi)\|^2_{\ca V} =&\sum_{p=1}^{m_a}\sum_{p'=1}^{m_a}\kappa_p(\xi)\kappa_{p'}(\xi)(\bar{{\ca L}}^p, \bar{{\ca L}}^{p'})_{\ca V}+\sum_{q=1}^{k-1}\eta_q(\xi)\sum_{p=1}^{m_a}\kappa_p(\xi)\\
&\times\bigg\{2\sum_{p'=1}^{m_a}\kappa_{p'}(\xi)(\bar{{\ca L}}^{p'}, {\ca L}_q^p)_{\ca V}
+\sum_{q'=1}^{k-1}\eta_{q'}(\xi)\sum_{p'=1}^{m_a}\kappa_{p'}(\xi)({\ca L}_q^{p}, {\ca L}_{q'}^{p'})_{\ca V}\bigg\}.
\end{split}
\end{equation}
and ${\ca L}_q^p$ is the Riesz representation of $a_p(\Phi_q, v)$, i.e., $({\ca L}_q^p, v)_{\ca V}=-a_p(\Phi_q, v)$ for any $v \in {\ca V}$ and $1 \leq p \leq m_a$, $1 \leq q \leq k-1$. Similarly, $\bar{{\ca L}}^p$ is the Riesz representation of $a_p({\bar\phi}, v)$, for any $v \in {\ca V}$ and $1 \leq p \leq m_a$. The derivation and calculation of $\|\hat e(\xi)\|^2_{\ca V}$ is similar to the one in paper \cite{li2017novel}.

Let $e(\xi)=e_h(x)e_\xi(\xi)$, and $e_h(x)$ be the solution of Eq. $(\ref{e-solution})$ with $\xi=\xi^k$.  If we  take $\Phi_k(x)=e_h(x)$ in Eq. $(\ref{vs-form})$, then
\begin{equation}
\label{e-form}
e_\xi(\xi)\sum_{p=1}^{m_a} \kappa_p(\xi)a_p(e_h, v)=r(v; \xi).
\end{equation}
If we take $v=e_h$ in Eq. $(\ref{e-form})$, then it follows that
\begin{equation}
\label{e-xi}
e_\xi(\xi)=\frac{-\sum_{p=1}^{m_a}\kappa_p(\xi)a_p({\bar\phi}, e_h)-\sum_{q=1}^{k-1}\eta_q(\xi)\sum_{p=1}^{m_a}a_p(\Phi_q, e_h)}{\sum_{p=1}^{m_a}\kappa_p(\xi)a_p(e_h, e_h)}.
\end{equation}
 We take $\eta_k(\xi)=e_\xi(\xi)$ in Eq.$(\ref{vs-form})$. As can be seen from Eq.$(\ref{e-xi})$, $\eta_k(\xi)$ depends on $\{\eta_q(\xi)\}_{q=1}^{k-1}$ computed previously while they are stochastic functions with respect to $\xi$. Let ${\mathbf b}^k:=[a_{1}({\bar\phi}, \phi^0_{k}), \cdots, a_{m_a}({\bar\phi}, \phi^0_{k})]$, $\kappa_\xi({\xi}):=[\kappa_{1}({\xi}), \cdots, \kappa_{m_{a}}({\xi})]^{T}$, and $\mathbf{\eta}^k({\xi}):=[\eta_1({\xi}), \cdots, \eta_{k-1}({\xi})]$,
\begin{eqnarray*}
{\mathbf A}^k:=\left[ \begin{array}{ccc}
 a_{1}(\phi^0_{1}, \phi^0_{1})& \ldots &  a_{m_a}(\phi^0_{1}, \phi^0_{k})\\ \vdots &\ddots &\vdots\\ a_{1}(\phi^0_{k}, \phi^0_{k})& \ldots & a_{m_a}(\phi^0_{k}, \phi^0_{k})\\\end{array}\right].
\end{eqnarray*}
Then the matrix form of $\eta_k({\xi})$ can be written as
\begin{eqnarray*}
\eta_k({\xi}) =\frac{{\mathbf b}^k\kappa_\xi(\xi)-{\eta}^k({\xi}){\mathbf A}_1^k\kappa_\xi(\xi)} {{\mathbf {\bar a}}^k\kappa_\xi(\xi)},
\end{eqnarray*}
where ${\mathbf A}_1^k$ is the first $k-1$ rows of matrix ${\mathbf A}^k$, ${\mathbf {\bar a}}^k$ is the last row of matrix ${\mathbf A}^k$. Because all the inner product in physical space can be computed and saved at each step, the simulation for $\xi$ is efficient due to the affine form.

\subsection{Ensemble-based VS method}

As can be seen, multiple inhomogeneous Dirichlet boundary condition problems are required to determined on each coarse neighborhood $\omega_i$. It will demand significantly computational cost if we approximate each solution $\phi^i_j(x, \xi)$ for $j=1, \cdots, M_i$ applying the VS method. Inspired by the idea of ensemble-based method used in  \cite{gunzburger2017ensemble}, we construct variable separated forms of $\phi^i_j(x, \xi)$ for different $j$ on each coarse neighborhood, such that they can share the same interpolate rule $\{\eta^{\omega_i}_q(\xi)\}_{q=1}^N$ .

By summing  up the Eq. $(\ref{effMS-basis})$ for different $j$, we can obtain the local problem
\begin{eqnarray}
\label{ensMS-basis}
\begin{cases}
& -\text{div}(\kappa(x; \xi)\nabla \langle\phi^i\rangle) =0, \ \ \text{in}\ \omega_i,\\
&\langle\phi^i\rangle = \langle\bar{\phi}^i(x)\rangle, \ \ \text{on}\ \partial \omega_i,
\end{cases}
\end{eqnarray}
where $\langle\bar{\phi}^i(x)\rangle=\sum_{j=1}^{M_i}\bar{\phi}^i_j$, and $\langle\phi^i(x)\rangle=\sum_{j=1}^{M_i}\phi^i_j$, $M_i$ is the number of local basis functions selected in the coarse neighborhood $\omega_i$. $\langle\phi^i(x)\rangle$ is treated as the solution to the new Eq. $(\ref{ensMS-basis})$, of which the inhomogeneous Dirichlet boundary condition is $\langle\bar{\phi}^i(x)\rangle$. Thus we have the solution
\[
\langle\phi^i(x)\rangle=\langle\phi^i(x)\rangle^0_N+\langle\bar{\phi}^i(x)\rangle,
\]
where $\langle\phi^i(x)\rangle^0_N$ has the variable separated approximation
\[
\langle\phi^i(x)\rangle^0_N=\sum_{q=1}^N\Phi^{\omega_i}_q(x)\eta^{\omega_i}_q(\xi).
\]
This can be  sought by running the VS procedure described in Section $\ref{VS}$. Here  the superscript 0 refers to the zero boundary condition, and the subscript $N$ represents the number of truncations terms.

At the step $k$ ($k\geq2$) of VS iteration, we solve the equation analogous to Eq. $(\ref{e-solution})$ with $\xi=\xi^k$, i.e.,
\[
a(\Phi^{\omega_i}_k(x), v; \xi^k)=-a(\langle\bar{\phi}^i(x)\rangle, v; \xi^k)-\sum_{q=1}^{k-1}\eta^{\omega_i}_q(\xi^k)a(\Phi^{\omega_i}_q, v), \ \forall  v\in {\ca V}_h,
\]
to obtain the $k-$th physical basis $\Phi^{\omega_i}_k(x)$. For $q=1, \cdots, k$, let
\[
\Phi^{\omega_i}_q(x)=\sum_{j=1}^{M_i}\Phi^{\omega_i}_{j,q}(x).
\]
Then the equation above can be rewritten as
\begin{equation}
\label{decom}
\sum_{j=1}^{M_i}a(\Phi^{\omega_i}_{j,k}(x), v; \xi^k)=-\sum_{j=1}^{M_i}a(\bar{\phi}^i_j(x), v; \xi^k)
-\sum_{q=1}^{k-1}\eta^{\omega_i}_q(\xi^k)\sum_{j=1}^{M_i}a(\Phi^{\omega_i}_{j,q}(x), v), \ \forall  v\in {\ca V}_h.
\end{equation}
We decompose Eq. $(\ref{decom})$ to be the system
\begin{eqnarray}
\label{decom-s}
\begin{split}
a(\Phi^{\omega_i}_{1,k}(x), v; \xi^k)&=-a(\bar{\phi}^i_1(x), v; \xi^k)-\sum_{q=1}^{k-1}\eta^{\omega_i}_q(\xi^k)a(\Phi^{\omega_i}_{1,q}(x), v),\\
\vdots& \\
a(\Phi^{\omega_i}_{M_i,k}(x), v; \xi^k) &=-a(\bar{\phi}^i_{M_i}(x), v; \xi^k)-\sum_{q=1}^{k-1}\eta^{\omega_i}_q(\xi^k)a(\Phi^{\omega_i}_{M_i,q}(x), v).
\end{split}
\end{eqnarray}
We note that  Eq. $(\ref{decom-s})$ is one of the decompositions. It seems like the residual equation of problem $(\ref{effMS-basis})$ for different $j$
in the $k-$th step of VS method, but the samples $\{\xi^q\}_{q=1}^k$ and interpolate rule $\{\eta^{\omega_i}_q(\xi)\}_{q=1}^{k-1}$ has been determined based  on problem $(\ref{ensMS-basis})$.
Let $\Phi^{\omega_i}_{j,k}(x)$ be the $k-$th physical basis of the corresponding variable separated form of the $j-$th local problem. We then obtain the separated form for each $\phi^i_j(x; \xi)$
\begin{equation}
\label{Enform}
\phi^i_j(x; \xi)=\sum_{q=1}^N\Phi^{\omega_i}_{j,q}(x)\eta^{\omega_i}_q(\xi), \ \ \text{for}\ \ 1\leq j\leq M_i,\ \ x\in \omega_i,
\end{equation}
which is considered as approximations of $\phi^i_j(x; \xi)$ satisfying Eq. $(\ref{effMS-basis})$ for $x\in \omega_i$.
\begin{rem}
The ensemble-based method we used here leads to a rough approximation of solutions to $(\ref{effMS-basis})$. We use the VS method to approximate the local problems with average boundary condition $\langle\bar{\phi}^i(x)\rangle$, and then derive the separated forms (mainly the physical parts) for each $\phi^i_j(x; \xi)$ through decomposing the residual equation.
\end{rem}
\begin{rem}
The approximation may not be as good as running VS for each $\phi^i_j(x; \xi)$, but the compromised scheme makes them share the same interpolate rule. The derived expression $(\ref{Enform})$ is an approximation of $\phi^i_j(x; \xi)$ in some sense. As $\Phi^{\omega_i}_{j,k}(x)$ are fixed, the  use of a common  interpolate rule significantly improves the computation efficiency.
\end{rem}
\begin{rem}
We can obtain the physical basis functions for each $j$ along with the interpolate rule $\{\eta^{\omega_i}_q(\xi)\}_{q=1}^N$. Or as an alternative, we run VS process for problem $(\ref{ensMS-basis})$ to obtain the selected samples $\{\xi^q\}_{q=1}^N$ and interpolate rule $\{\eta^{\omega_i}_q(\xi)\}_{q=1}^N$. Then calculate the physical basis by the similar technique used  in paper \cite{ou2019new}, for any $j=1, \cdots, M_i$.
\end{rem}

We can finally multiply the physical basis functions with the  partition functions $\{\chi_i(x)\}_{i=1}^{N_H}$, and then the stochastic multiscale basis functions are presented as
\[
\phi^i_j(x; \xi)=\sum_{q=1}^N\chi_i(x)\Phi^{\omega_i}_{j,q}(x)\eta^{\omega_i}_q(\xi), \ \ \text{for}\ \ 1\leq j\leq M_i, \ \ x\in \omega_i.
\]
In particular, we set $\phi^i_j(x; \xi)={\chi}_i(x)\bar{\phi}^i_j(x)$ on $\partial \omega_i$. The physical part $\chi_i(x)\Phi^{\omega_i}_{j,k}(x)$ can be calculated and stored for any random sample $\xi$. Thus  the stochastic multiscale basis functions can be efficiently obtained. Algorithm $\ref{EnGMsVS}$ describes the procedure of ensemble-based VS-Ms method.

\begin{algorithm}
  \caption{Ensemble-based VS-Ms method}
1. Take an arbitrary realization $\log{\kappa(x; \bar{\xi})}$ from the prior $(\ref{GP})$, and generate the effective boundary conditions $\bar{\phi}_j^i(x)$ for $j=1, \cdots, M_i$ by GMsFEM in coarse neighborhoods;\\
2. Run Ensemble-based VS procedure for local problems $(\ref{ensMS-basis})$, to obtain the interpolate rules $\{\eta^{\omega_i}_q(\xi)\}_{q=1}^N$;\\
3. In each coarse block, reverse the physic basis for $\phi^i_j(x; \xi)$ and obtain the corresponding separated form \\
\[
\phi^i_j(x; \xi)=
\begin{cases}
\sum_{q=1}^N\Phi^{\omega_i}_{j,q}(x)\eta^{\omega_i}_q(\xi), \ \ & x\in \omega_i,\\
\bar{\phi}^i_j(x), \ & x\in \partial \omega_i.
 \end{cases}
\]

4. Multiply the physical basis with the union partition functions, and obtain the stochastic multiscale basis functions as \\
\[
\phi^i_j(x; \xi)=
\begin{cases}
\sum_{q=1}^N\chi_i(x)\Phi^{\omega_i}_{j,q}(x)\eta^{\omega_i}_q(\xi), \ \ &x\in \omega_i,\\
{\chi}_i(x)\bar{\phi}^i_j(x), \ & x\in \partial \omega_i.
\end{cases}
\]
5. For any random realization taken from the prior $(\ref{GP})$, calculate the multiscale basis functions efficiently and simulate the model.
    \label{EnGMsVS}
\end{algorithm}

The random samples of   $\xi$ are taken from the prior $(\ref{GP})$, where the covariance matrix is determined by the relative distance of the spatial points. Due to the property of the Gaussian process, the parts of $\xi$ taken from the same local domain structures obey the same Gaussian distribution. In this paper, $\bar \xi$ is chosen as the  mean of the random variable $\xi$,
the local boundary conditions $\bar{\phi}^i_j(x)$ for $j=1, \cdots, M_i$ are  the same for the same types of local structures. Hence, we only need to use the proposed method to obtain stochastic multiscale basis functions for 4 types of local domain structures (as shown in red box of Figure $\ref{coarse-cell}$). We apply the corresponding 4 types of expressions to the whole domain $D$.

\section{Numerical examples}

In this section, we consider the discontinuous field identification problems for a single-phase Darcy flow problem. The flow is described in terms of fluid pressure $u(x, t)$ which is governed by the following equation
\begin{equation}\label{PPDE}
c\frac{\partial u}{\partial t}-\text{div}\bigg(\kappa(x)\nabla u\bigg)=f(x, t),\ \ x\in D, t\in(0,T],
\end{equation}
where the studied domain is $D= [0, 1] \times [0, 1]$, and the initial condition is set as $u(x,0)=0$. $\kappa(x)$ is the unknown permeability field. $c$ is the specific storage, which is taken to be constant 1 for simplicity.
The source term $f(x, t)$ is defined by a weighted Gaussian plume with standard deviation 0.1, centered at $(0.55, 0.4)$, with weight 1, which keeps constant with respect to the time $t$. The log permeability is parameterized as a Gaussian random field, the parameters in the Gaussian kernel $\ca{C}(x_1, x_2)$ are set as
\[
\sigma_x^2=0.1, \ l_{x_1}=0.07,\ l_{x_2}=0.07.
\]

Eq. $(\ref{affinem})$ is treated as the affine form of the field $\log{\kappa(x; \xi)}$, to characterize the unknown random field accurately. Then we have the required affine of the bilinear form
\[
a(w, v; \xi) =\sum_{p=1}^{m_a}e^{\xi_p}a_p(w, v), \ \forall v,w \in {\ca V}, \ \forall \xi \sim {\mathcal N}(0, \Sigma),\\
\]
where
\[
a_p(w, v)=\int_D \mathbf{I}_{\{x\in D_i\}}\nabla w \nabla v dx.
\]
Note that there are only a fixed number of nonzeros in the corresponding stiffness matrix for each $a_p(w, v)$, we use the proposed method on 4 types of local coarse neighborhood, the dimension $m_a$ would not be too large to slow down the computation. For the simplicity  of notation, let $m(x; \xi)=\log{\kappa(x; \xi)}$. The backward Euler scheme is used for the temporal discretization. Some model simulations are presented in the first subsection to demonstrate the efficiency of the scheme $(\ref{effMS-basis})$ and our proposed method in solving the forward model. We also apply the proposed model reduction method to the  inverse problems to show the  performance of the method in identifying discontinuous fields.

\subsection{Reduced model approximation}
\label{reductionapp}

We study the problem $(\ref{PPDE})$ with no-flow boundary conditions. The end time is $T=0.2$, and the time step is set as $\Delta t=0.001$. The solution obtained by finite element method is considered the reference one. We solve it on a uniform $80\times80$ fine grid, and we set the coarse grid  $8\times8$ for multiscale finite element method.

The local boundary conditions in problem $(\ref{effMS-basis})$ are generated by setting $m(x; \bar{\xi})$ to be the mean value of the random field. The ensemble-based VS-Ms method is used to obtain the variable separated form of multiscale basis functions, where 4 local stochastic problems are solved on 4 coarse neighborhood as shown in Figure $(\ref{coarse-cell})$. The parameter dimension $m_a$ on the coarse neighborhood $\mathrm{I}$, $\mathrm{II}$, $\mathrm{III}$, $\mathrm{IV}$ is 100, 200, 200, and 400, respectively.

We show the effects of the ensemble-based VS-Ms method in solving the problem $(\ref{PPDE})$ numerically.
 Let $u_h$ be the reference solution using the  finite element method on fine grid, $u_H$ the solution obtained by GMsFEM. $u_{\bar m}$ represents the solution obtained by the multiscale method, where the reduced snapshot functions are obtained by solving local problems $(\ref{effMS-basis})$, $u_{\langle \bar m\rangle}$ is the approximation resulted from the ensemble-based VS-Ms method. Let $\varepsilon^n_s=u_h^n-u_s^n$, $s=H,{\bar m},{\langle \bar m\rangle}$. Then the $L^\infty$-error and $L_2$-error at $n\Delta t$ are defined as
\[
\|\varepsilon^n_s\|_\infty= \max_{\xi\in\Omega}\max_{x\in D}|u^n_h(x; \xi)-u^n_s(x; \xi)|,
\]
\[
\|\varepsilon^n_s\|=\sqrt{\mathlarger{\int}_\Omega\mathlarger{\int}_D (u^n_h(x; \xi)-u^n_s(x; \xi))^2dxd\pi_0}, \ s=H,{\bar m},{\langle \bar m\rangle},
\]
where $\pi_0$ is the Gaussian prior ${\mathcal N}(0, \Sigma)$.  We compare the maximum $L^\infty$-error and the maximum $L_2$-error for the different strategies, and define
\[
\varepsilon_{\infty, s}=\max_{0\leq n\leq n_T}\|\varepsilon^n_s\|_\infty,
\]
\[
\varepsilon_{2, s}=\max_{0\leq n\leq n_T}\|\varepsilon^n_s\|,
\]
to show the corresponding performance, where $T=n_T\Delta t$.

For $s=H,{\bar m},{\langle \bar m \rangle}$, $\tau_s$ represents the time cost in constructing the multiscale basis functions for an arbitrary realization of $\xi$. The numerical results are shown in Table $\ref{Error}$.  $81 $ local problems are required to be solved for all the methods displayed above. As for  GMsFEM, the effective boundary conditions are obtained by solving eigenvalue problems $(\ref{redsnapshots})$, which depends on  the stochastic random field. Solving local problems $(\ref{effMS-basis})$ directly or using the ensemble-based VS-Ms method cost much less time, since the local boundary conditions are fixed and can be stored overhead. The separated form provided by the ensemble-based VS-Ms method significantly  accelerates the speed of calculating multiscale basis functions.

\begin{table}
  \centering
  \begin{tabular}{|c|c|c|c|}
    \hline
    number of local basis functions  & 3 & 5 & 7   \\  \hline  
    $\tau_H$ & 0.4359 s& 0.4456 s& 0.4475 s\\
    ${\varepsilon }_{\infty, H}$ & 3.9581e-03 & 3.7005e-03 & 3.6726e-03 \\
    $\varepsilon_{2, H}$ & 1.4300e-04 & 8.7137e-05 & 4.1381e-05  \\ \cline{1-4}
    $\tau_{\bar m}$ & 0.0326 s& 0.0396 s& 0.0469 s\\
    $\varepsilon_{\infty, {\bar m}}$ & 3.7533e-03 & 3.5010e-03  & 3.4626e-03 \\
    $\varepsilon_{2, {\bar m}}$ & 1.7949e-04 & 1.2046e-04 &  5.7366e-05\\ \cline{1-4}
    $\tau_{\langle \bar m \rangle}$ & 0.0076 s&  0.0137 s& 0.0177 s\\
    $\varepsilon_{\infty, {\langle \bar m \rangle}}$ & 4.8471e-03  & 3.5276e-03 & 3.01960e-03 \\
    $\varepsilon_{2, {\langle \bar m \rangle}}$ & 2.5245e-04 & 1.4574e-04 & 8.3739e-05 \\
    \hline
  \end{tabular}
  \caption{The maximum $L_\infty$,  maximum $L_2$ error and time cost in constructing the local basis functions each time, with different number of selected basis functions in each coarse neighborhood,  the separated terms is fixed at 20. }
  \label{Error}
\end{table}

As can be seen from the table, both the two types of errors decrease as the number of multiscae basis functions increase, for the studied strategies. Under the prior condition, the scheme $(\ref{effMS-basis})$ we adopt in generating the reduced snapshot space works well, it reaches almost the same magnitude accuracy as GMsFEM. The accuracy of the PDE is preserved when using the ensemble-based VS-Ms method. Due to the constructed explicit expressions, the basis functions can be efficiently computed. Hence  it costs the least time in constructing local basis functions.

When the number of multiscale basis functions is fixed by $3$, we solve the problem $(\ref{PPDE})$ by varying $N=5$, 15, 25, 35, to study the effect of the number of separated terms. It can be shown from Table $\ref{ErrorN}$ that, the trend of the two types of errors is flat against the number of separated terms, which is resulted from the high efficiency of the VS method with small separated terms. The number of multiscale basis functions dominates the effects on the accuracy for our proposed method. Hence, we will concentrate on the influence of the number of multiscale basis functions on the posterior density approximation in the next section.

\begin{rem}
Though the problem can be solved well with even small separated terms, the moving of samplers in the Markov chain leads them away from the prior, we set $N$ not be too small to make the proposed method stable and still work in Markov chain sampling.
\end{rem}

\begin{table}
 \centering
 \begin{tabular}{|c|c|c|c|c|}
   \hline
   separated terms $N$& 5 & 15 & 25 & 35 \\ \hline
   $\tau_{\langle \bar m \rangle}$ &  0.0053&  0.0063& 0.0081 & 0.0106 \\
   $\varepsilon_{\infty, {\langle \bar m \rangle}}$ & 4.8438e-03 & 4.8378e-03 & 4.8471e-03 & 4.8410e-03\\
   $\varepsilon_{2, {\langle \bar m \rangle}}$ & 2.5241e-04 & 2.5267e-04 & 2.5257e-04 & 2.5271e-04 \\
   \hline
 \end{tabular}
 \caption{The maximum $L_\infty$,  maximum $L_2$ error and time cost in constructing the multiscale basis functions each time, with different number of separated terms, the number of multiscale basis functions is fixed at 3. }
 \label{ErrorN}
\end{table}


\subsection{Posterior densities approximation}

In this section, we use the proposed method to approximate the likelihood, along with the imposed TG hybrid prior to identify  discontinuous fields. We construct the variable separated expression for the multiscale basis functions by the ensemble-based VS-Ms method, and then use it to improve the efficiency of MCMC.
 Measurement data are generated by using finite element method in a fine grid with time step $\Delta t = 0.001$. For any given values of $\xi$, we solve the equation using the ensemble-based VS-Ms model reduction method with time step $\Delta t=0.002$. The measurement noise is set to be $\sigma=0.01$.

\subsubsection{Case 1: separated blocks}
\label{sepfacs}\begin{figure}[t]
  \centering
  \subfigure[]{
  \includegraphics[width=0.8\textwidth]{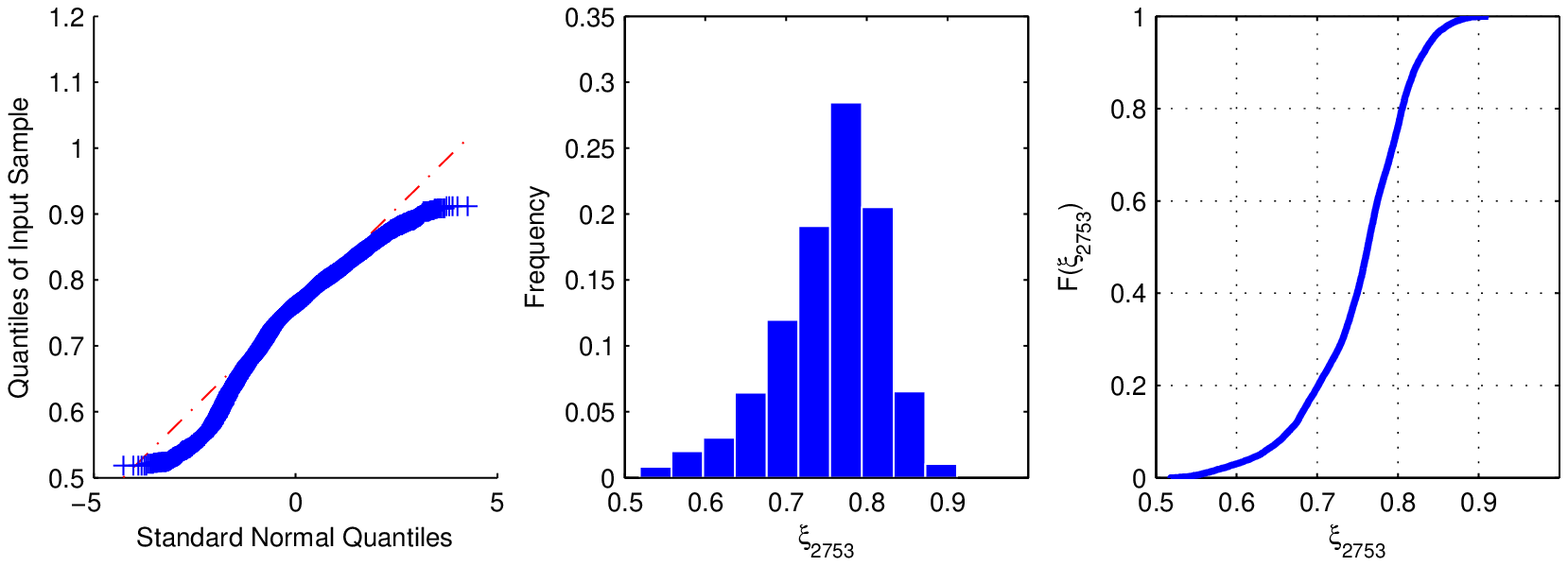}}
  \subfigure[]{
  \includegraphics[width=0.8\textwidth]{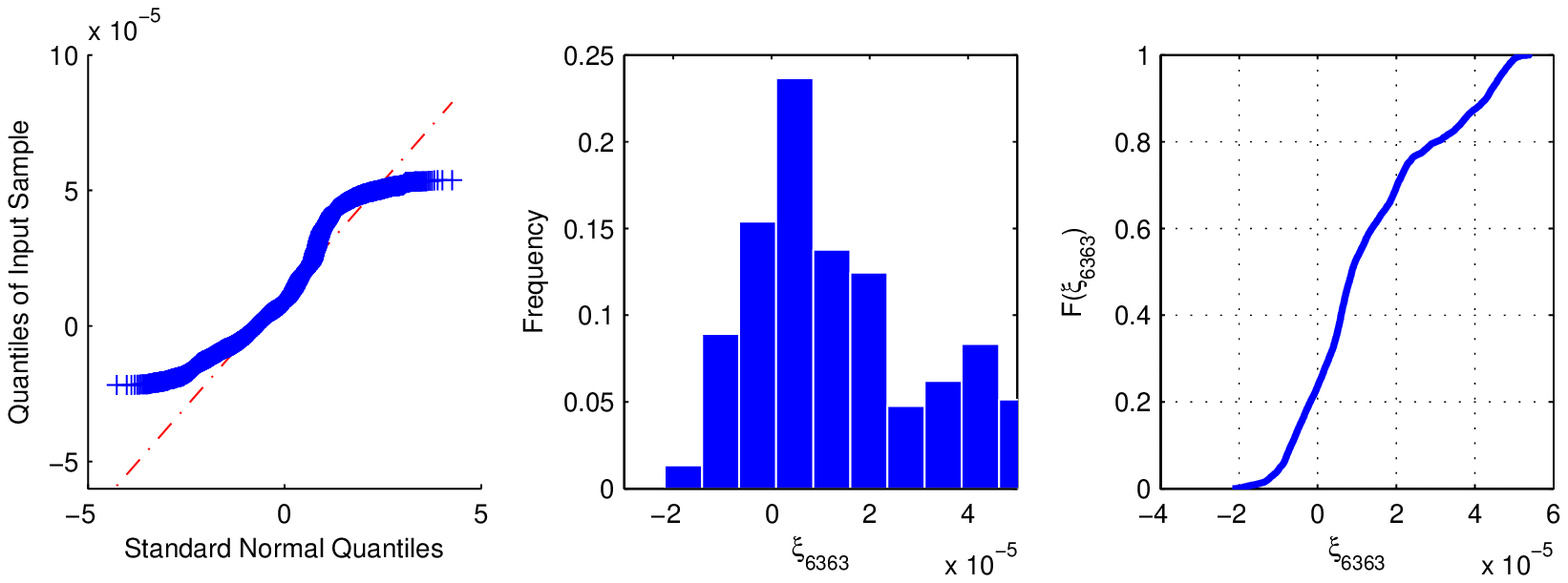}}
    \subfigure[]{
  \includegraphics[width=0.8\textwidth]{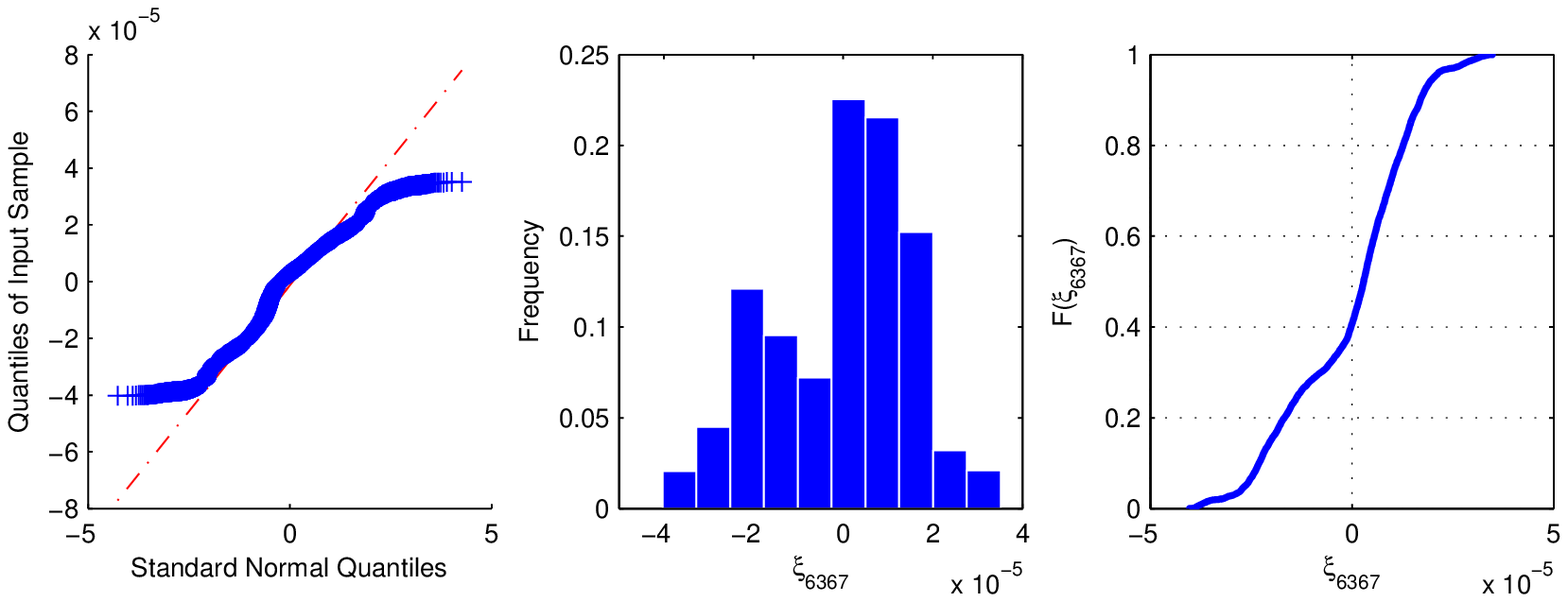}}
    \subfigure[]{
  \includegraphics[width=0.8\textwidth]{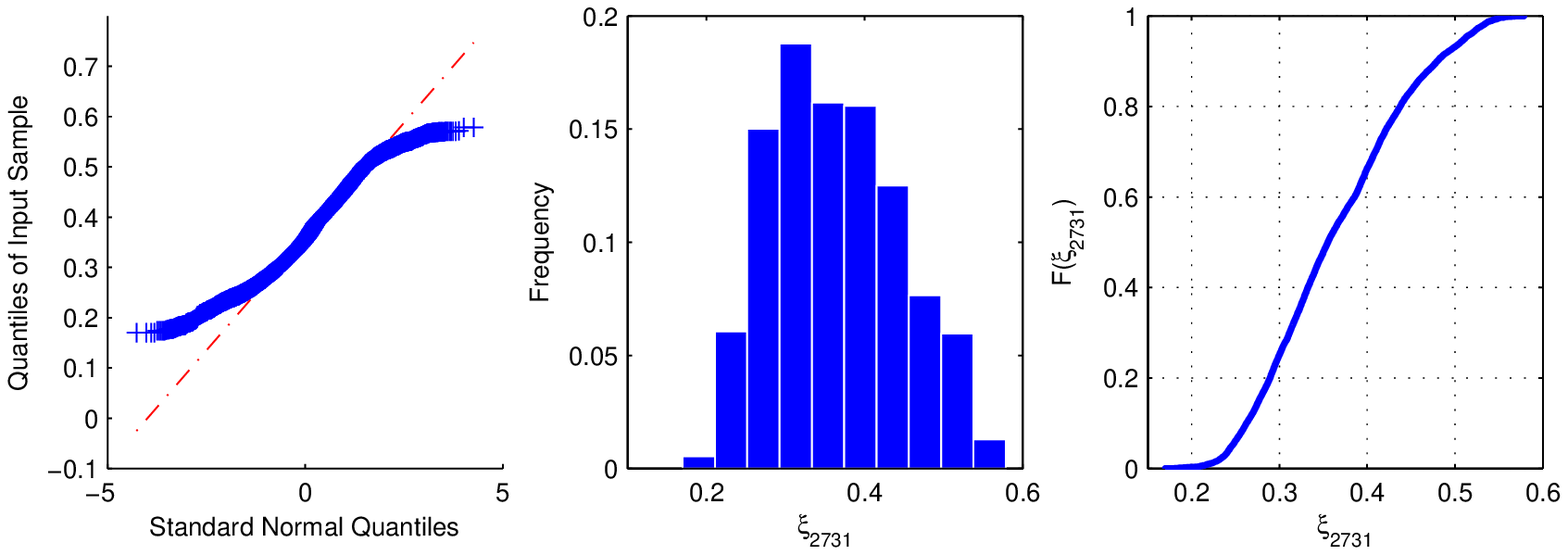}}
  \caption{The empirical quantile-quantile, histogram and empirical cdf of (a) $\xi_{2753}$, (b) $\xi_{6363}$, (c) $\xi_{6367}$ and (d) $\xi_{2731}$ .}
  \label{qqplot}
\end{figure}

First we apply the model reduction method to the  example, where there are two separated blocks on the target permeability field. Namely, suppose we have the target piecewise constant permeability
\[
\kappa(x)=\begin{cases}
e^2, & \frac{(x_1-0.35)^2}{0.15^2}+\frac{(x_2-0.3)^2}{0.15^2}\leq 1,\\
e^2, &\frac{(x_1-0.65)^2}{0.15^2}+\frac{(x_2-0.6)^2}{0.15^2}\leq 1,\\
1, & \text{otherwise}.
\end{cases}
\]
i.e., there are two faces to be identified. The boundary condition is set as
\[
u(x, t)=1.7-1.4x_1, \ \ x \in \partial D,
\]
for the problem, the forward model is studied during the time interval $[0, 0.15]$. We observe the data starting from $t=0.02$, at every 0.01, up to $t=0.11$. For each time layer, 49 pressure measurements are uniformly collected in $[0.05, 0.95]\times[0.05, 0.95]$.

The governing equation in $(\ref{PPDE})$ is solved on a uniform $80\times80$ fine grid, resulting in 6400 unknowns to be estimated. The coarse grid is set to be $8\times8$, the effective local boundary conditions are generated by setting $m(x; \bar{\xi})=0$, we select 3 basis functions on each coarse neighborhood, 20 separated terms are kept in the variable separation expressions. The posterior samplers are drawn with $\lambda=300$, $\beta=0.015$, resulting in the acceptance rate 38.55\%. We draw $10^5$ samplers to ensure the reliability of the inference, and remain the last 50\% for statistics computation.

We display the empirical quantile-quantile, histogram and empirical cumulative distribution function (cdf) of posterior samplers to show the non-Gaussian of the posterior density. The plots of $\xi_{2753}$, $\xi_{6363}$, $\xi_{6367}$ and $\xi_{2731}$ are shown in Figure $\ref{qqplot}$. The quantiles illustrate that the distributions differ from Gaussian significantly. As can be seen from the histogram, the distributions of $\xi_{2753}$ and $\xi_{2731}$ are right-skewed and left-skewed, respectively. For $\xi_{2731}$, the cdf curve reaches 0.5 fast, but the interval gets wider when it reaches  the rest 0.5. $\xi_{6363}$ and $\xi_{6367}$ seem to be multi-modal, which obviously are non-Gaussian.

The true profile, posterior mean and posterior standard deviation of the log permeability field are shown in Figure $\ref{block2}$. The two blocks are separately identified, no adjoin between them exists as expected.
Both the values of the two blocks are underestimated, the block near the lower left corner is estimated better. It may be resulted by the higher values condition imposed at the left boundary when solving the forward problem.
\begin{figure}
  \centering
  \includegraphics[width=0.9\textwidth,height=0.3\textwidth]{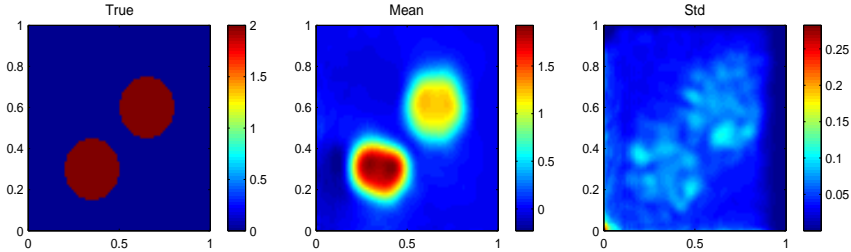}\\
  \caption{True profile, posterior mean and posterior standard deviation of the log permeability field.}
  \label{block2}
\end{figure}

\begin{figure}
  \centering
  \subfigure[]{
  \includegraphics[width=0.6\textwidth]{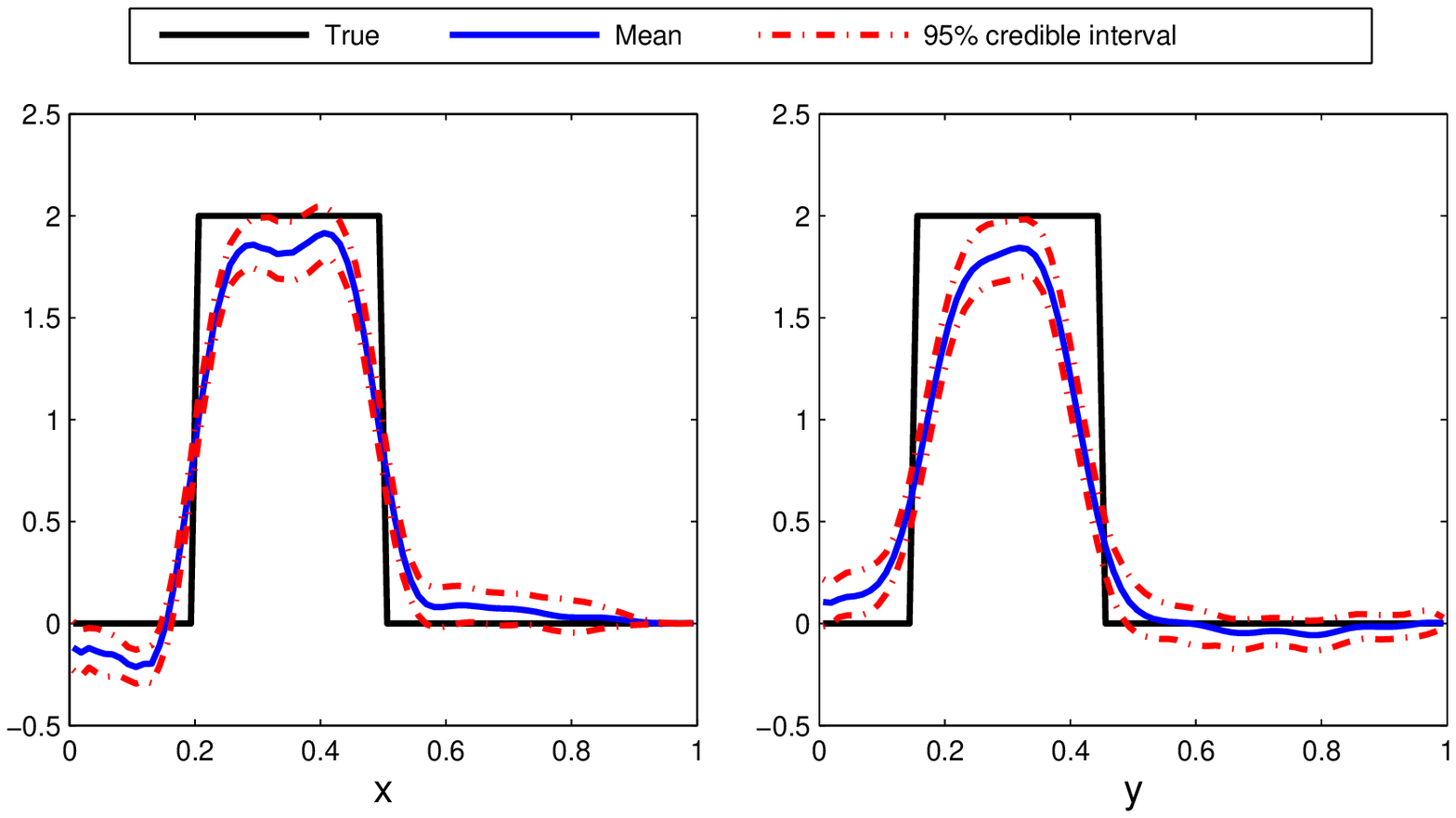}}
  \subfigure[]{
  \includegraphics[width=0.6\textwidth]{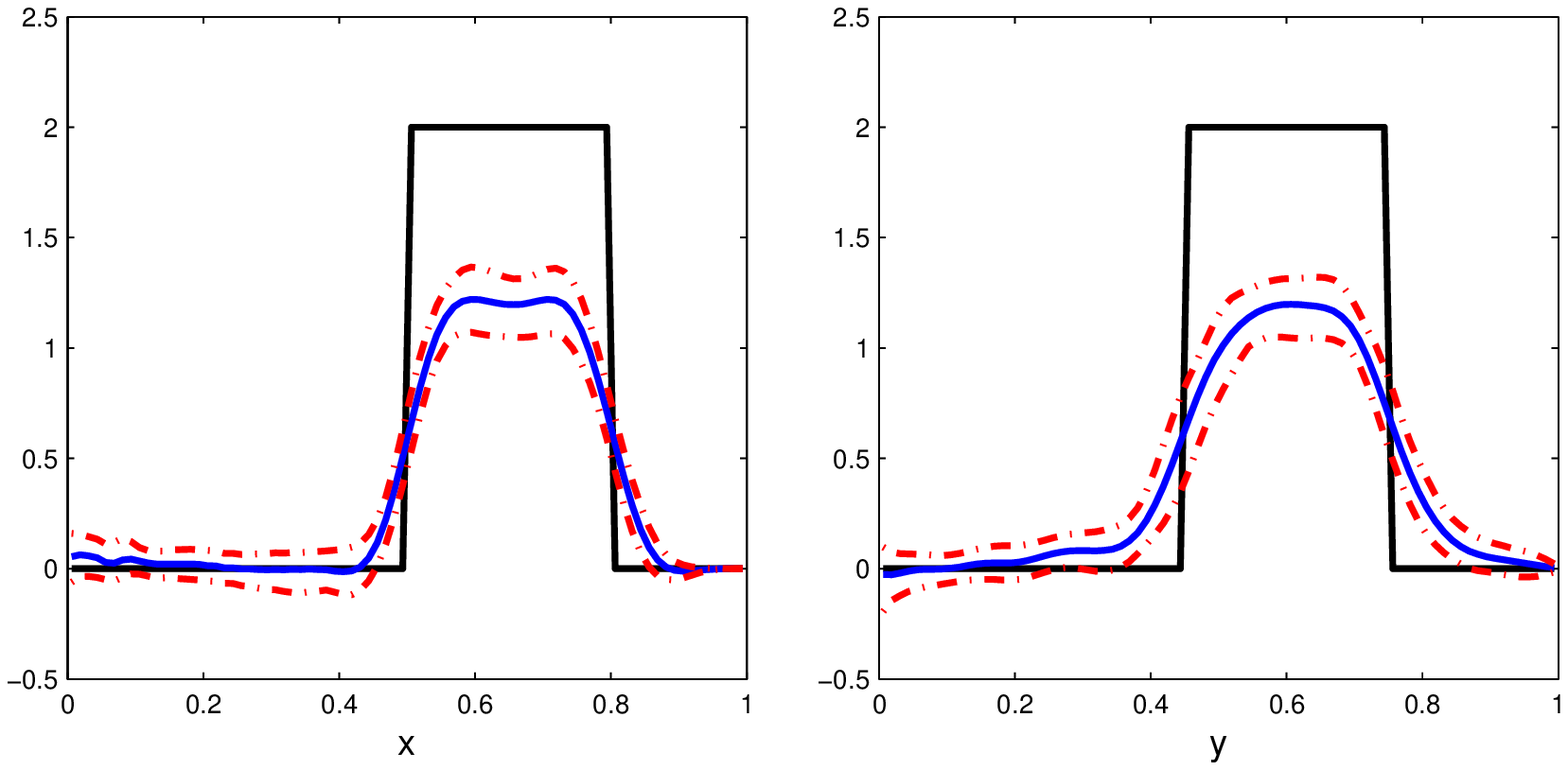}}
  \caption{The true profile, posterior mean and 95\% credible estimates of the log permeability at the surface (a) $y=0.3$ (left) and $x=0.35$ (right), (b) $y=0.6$ (left) and $x=0.65$ (right).}
  \label{credible2}
\end{figure}

\begin{figure}
  \centering
  \includegraphics[width=0.6\textwidth]{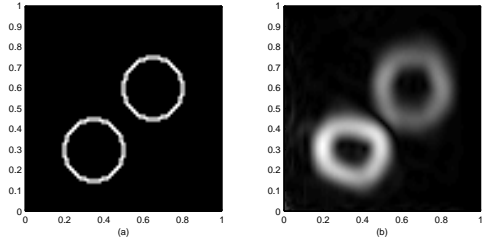}\\
  \caption{Spatial gradient of (a) true and (b) posterior mean of the log permeability field.}
  \label{gradient2}
\end{figure}

\begin{figure}
  \centering
  \includegraphics[width=1\textwidth]{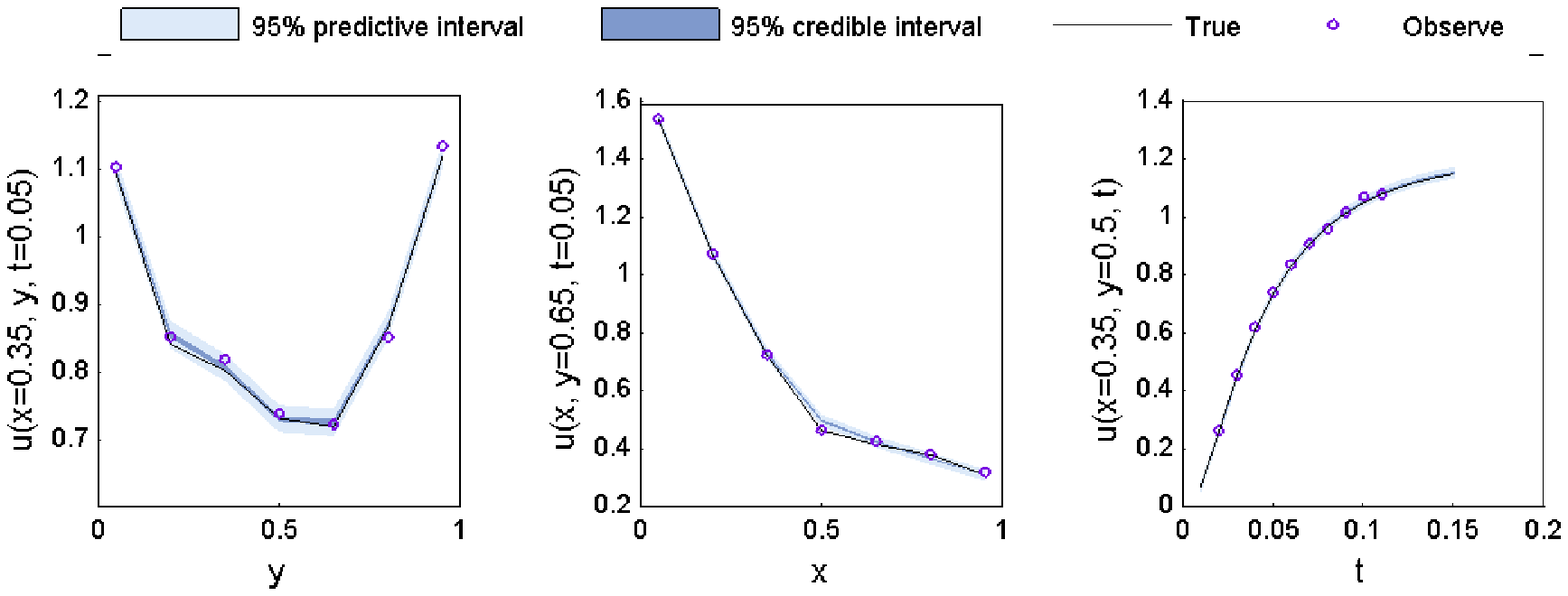}\\
  \caption{Data, point estimates, and 95\% credible and prediction intervals produced by the Bayesian analysis for $u\big(x=0.35, y, t=0.05\big)$,
 $u\big(x, y=0.65, t=0.05\big)$ and $u\big(x=0.35, y=0.5, t\big)$.}
  \label{pre2}
\end{figure}

\begin{figure}
  \centering
  \includegraphics[width=0.6\textwidth]{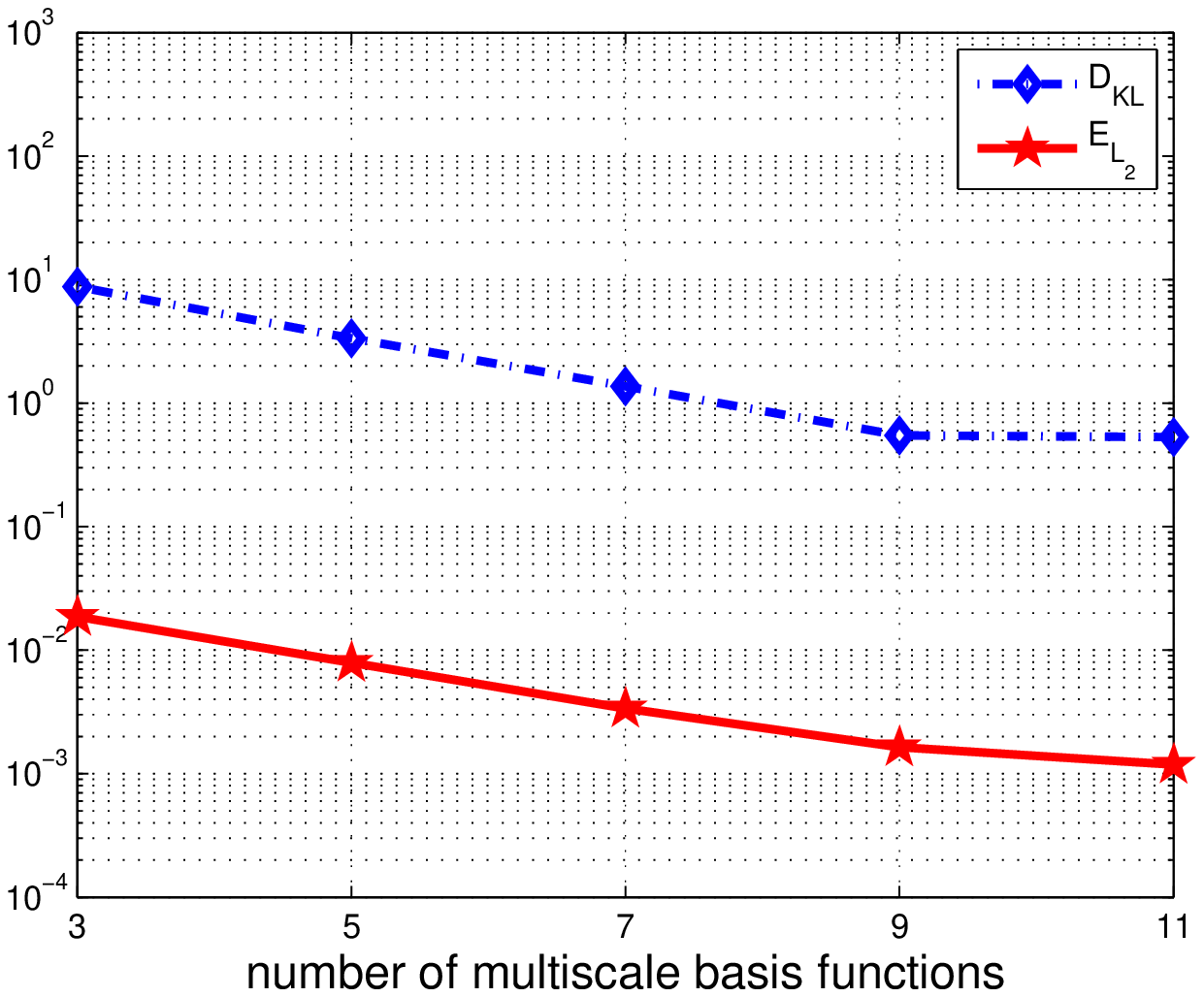}\\
  \caption{Approximation of the forward model and the posterior density with respect to the number of multiscale basis functions. Dashed line represents the Kullback-Leibler divergence $D_{KL}(\pi_{\bff N}^{d}||\pi^d)$ (denoted by $D_{KL}$); solid line shows the $L_{\pi_0}^2$ error  $\|{\bff G}-\bff {G_N}\|^2_{L_{\pi_0}^2}$ (denoted by $E_{L_2}$). }
  \label{KLD2}
\end{figure}
We plot the surface from the view of plane $y=0.3$, $y=0.6$, $x=0.35$, and $x=0.65$ in Figure $\ref{credible2}$, the mean and 95\% credible intervals are presented as well.  The shapes of the curves imply that the background and jumps are recognized well, and  the value of the log permeability is very closed to the true one for the block near the lower left corner. Though the jumps estimated are not as sharp as the reference ones, there are no oscillations near them. It can also be reflected from Figure $\ref{gradient2}$, where the spatial gradient of the true and posterior mean of the log permeability field are presented. Compared with the reference face boundaries, the estimated ones are wider because of the less sharp estimated jumps. The face boundary near the upper right corner is darker since the severe underestimated values at this block. However, the locations of the face boundaries are separately identified generally.

The 95\% credible intervals and prediction intervals for the pressure at $u\big(x=0.35, y, t=0.05\big)$,
$u\big(x, y=0.65, t=0.05\big)$ and $u\big(x=0.35, y=0.5, t\big)$ are constructed. As illustrated in Figure $\ref{pre2}$, the reference values are mostly located in the credible intervals, measurements are almost contained in the predictive intervals. Due to the deterministic Dirichlet boundary condition on $\partial D$, both the credible and prediction intervals become tight as $y$ and $x$ gets closer to the end points for $u\big(x=0.35, y, t=0.05\big)$ and $u\big(x, y=0.65, t=0.05\big)$. For $u\big(x=0.35, y=0.5, t\big)$, as the uncertainty from the input propagates, the uncertainty associated with the model fit and predictions grows, which naturally leads to loose intervals as time moves on.

Figure $\ref{KLD2}$ shows convergence of the forward model and posterior approximation with respect to the number of multiscale basis functions selected in the coarse neighborhood. We plot the $L_2$ error in the forward model, $\|{\bff G}-\bff {G_N}\|^2_{L_{\pi_0}^2}$ and the Kullback-Leibler divergence of the exact posterior from the approximate posterior, $D_{KL}(\pi_{
\bff N}^{d}||\pi^d)$. Both the errors $D_{KL}$ and $E_{L_2}$ decrease exponentially as the number of multiscale basis functions increase. The convergence rate of the Kullback-Leibler divergence is basically the same as that of the $E_{L_2}$ error in the forward model. This confirms  Theorem $\ref{DKL}$.

\subsubsection{Case 2: nested blocks}

As the last example, we consider recovering a discontinuous permeability with nested blocks structure. The problem is studied with boundary condition
\[
u(x, t)=2-2x_1, \ \ x \in \partial D,
\]
and the end time $T=0.1$. The reference log permeability field and measurement locations are shown as Figure $\ref{RefLoc}$, where  the two blocks are nested with different permeability value. We observe data in the time interval $[0.02, 0.1]$ with the time step $0.01$.

\begin{figure}
  \centering
  \includegraphics[width=0.5\textwidth]{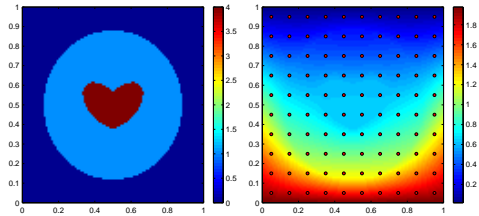}\\
  \caption{(Left) The reference log permeability field and (Right) true pressure field at time $t=0.05$, and the observation locations. }
  \label{RefLoc}
\end{figure}

\begin{figure}
  \centering
  \includegraphics[width=0.9\textwidth,height=0.3\textwidth]{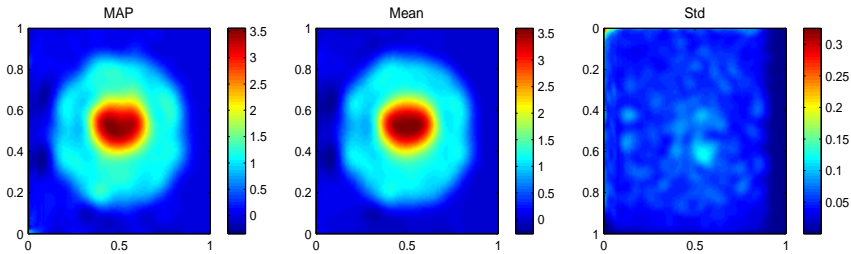}\\
  \caption{The MAP estimate, posterior mean and posterior standard deviation of the log permeability field.}
  \label{block3}
\end{figure}

\begin{figure}
  \centering
  \includegraphics[width=0.6\textwidth]{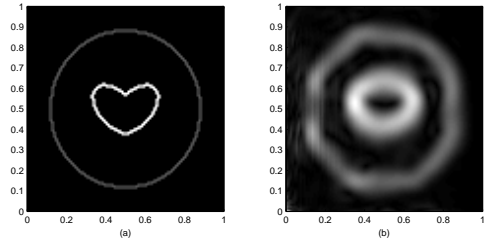}\\
  \caption{Spatial gradient of (a) true and (b) mean of the reconstructed log permeability field.}
  \label{gradient3}
\end{figure}
\begin{figure}
  \centering
  \includegraphics[width=0.6\textwidth]{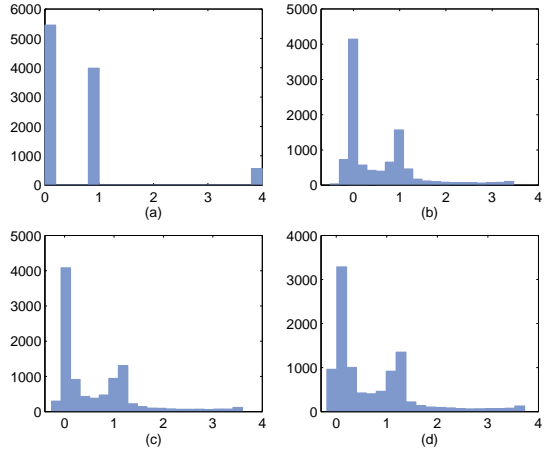}\\
  \caption{The histogram of the (a) reference, (b) 0.25 quantile, (c) mean and (d) 0.95 quantile estimate of the log permeability. }
  \label{hist}
\end{figure}

In this subsection, we solve the problem on a uniform $100\times100$ fine grid, which results in $10^4$ unknowns. For the sake of using the explicit variable separation expression of the multiscale basis functions constructed before, we set the coarse grid as $10\times10$. Then the structure of local covariance is the same, but the relative distance becomes smaller. Correspondingly, we put the influence of the change into the length scale in the covariance kernel, and modify the length scale to be $0.8\times l_{x_1}$ and $0.8\times l_{x_2}$ in this example. This will result in the same Gaussian prior as the two examples presented above. It is convenient to use the constructed variable separated expression but add some basis index and rearrange the corresponding basis functions.

The Markov chain is run by the proposed method, where the number of basis functions is 3. We draw $1.4\times 10^5$ posterior samplers with $\lambda=250$, $\beta=0.012$, and keep the last $5\times 10^4$ ones for statistical inference, the resulted acceptance rate is 38.80\%.
Figure $\ref{block3}$ displays the maximum a posterior (MAP) estimate, posterior mean and posterior standard deviation of the log permeability field. The uncertainty is small overall, the MAP reflects the concave part of the heart shaped block better than the posterior mean, the lower half of the heart shaped block are well embodied by both the two estimates. Besides the background, the circle block is also well identified, whose value is between the background and the heart shaped block.

We present the spatial gradient of the true and posterior mean of the log permeability field in Figure $\ref{gradient3}$.  The face boundary for the nested block is brighter than the circle's. This implies the spatial gradient value is larger, which matches the spatial gradient of the true log permeability. The histogram of the reference, 0.25 quantile, mean and 0.95 quantile estimate of the log permeability are demonstrated in Figure $\ref{hist}$. The reference has a trinomial distribution in the histogram, the high frequency value concentrates around 0 and 1 for the estimates, the stack around 4 or 3.5 is less prominent due to the underestimation of the permeability value.

In Figure $\ref{pre3}$, we plot the 95\% credible and predictive intervals for the model response at $u\big(x=0.35, y, t=0.05\big)$, $u\big(x, y=0.65, t=0.05\big)$ and $u\big(x=0.35, y=0.5, t\big)$. Most part of the reference values lie in the credible intervals, apart from one or two points, measurements are located in the predictive intervals. Similar to the response studied in Section $\ref{sepfacs}$, for $u\big(x=0.35, y, t=0.05\big)$ and $u\big(x, y=0.65, t=0.05\big)$, the uncertainty of the model response mainly concentrates in the domain, and gets  smaller as the response distributes closer to $\partial D$. The accumulated propagation uncertainty results in loose intervals for $u\big(x=0.35, y=0.5, t\big)$.

\begin{figure}
  \centering
  \includegraphics[width=1\textwidth]{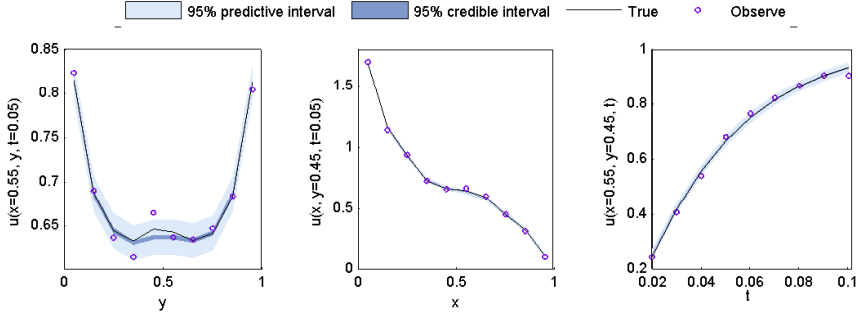}\\
  \caption{Data, point estimates, and 95\% credible and prediction intervals produced by the Bayesian analysis for $u\big(x=0.55, y, t=0.05\big)$,
 $u\big(x, y=0.45, t=0.1\big)$ and $u\big(x=0.55, y=0.45, t\big)$.}
  \label{pre3}
\end{figure}

\begin{figure}
  \centering
  \includegraphics[width=0.6\textwidth]{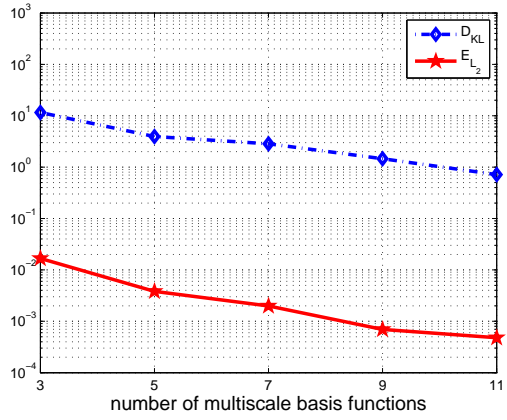}\\
  \caption{Approximation of the forward model and the posterior density with respect to the number of the multiscale basis functions. Dashed line represents the Kullback-Leibler divergence $D_{KL}(\pi_{\bff N}^{d}||\pi^d)$ (denoted by $D_{KL}$); solid line shows the $L_{\pi_0}^2$ error  $\|{\bff G}-\bff {G_N}\|^2_{L_{\pi_0}^2}$ (denoted by $E_{L_2}$). }
  \label{KLD3}
\end{figure}

The Kullback-Leibler divergence and the forward model $L_2$ norm error are plotted against the number of multiscale basis functions in Figure $\ref{KLD3}$. We observe the same trend for all curves, the errors decay as the number of multiscale basis functions increase. Moreover, the curve of the figure confirms that the posterior density of Kullback-Leibler divergence $D_{KL}$ decreases with the same speed as $\|{\bff G}-\bff {G_N}\|^2_{L_{\pi_0}^2}$, and the error $D_{KL}$ also can be bounded in terms of
the error $E_{L_2}$ in the forward model.


\section{Conclusion}

This paper presented the ensemble-based VS-Ms model reduction method with application in Bayesian inverse problems for discontinuous field identification. The stochastic multiscale basis functions are approximated by the ensemble-based VS method, with respect to the reference Gaussian measure. The adopted ensemble scheme and residual decomposition strategy lead to a common interpolate rule for local functions on the same coarse neighborhood, which improves the efficiency of basis calculation. The way of generating local problems and the relative distance characterized covariance structure simplify the basis construction problem to 4 local problems, coarse neighborhood with the same partition of grid share the same variable separated form. With the same local partitions, only some basis index and rearrange the corresponding basis functions are needed for refinement of the fine grid partition. The influence of changing the relative distance is considered as length scale modification. The presented approach leads to an accurate approximation of the full forward model while save much computation cost, due to the efficient multiscale basis functions construction and the degree of freedom reduction. Rigorous analysis is carried out for the approximation of model reduced
method in the Bayesian inverse problem with the TG prior. The numerical examples confirm that
the approximated posterior matches the reference one better and better as the enrichment of multiscale basis functions. We demonstrate our approach in the single-phase inversion problems, where the discontinuous permeability with different structures are identified. The constructed explicit expression of the basis functions can be applied to problems with different right-hand conditions, which is inherited from the GMsFEM, while the computational time is less than the GMsFEM.

\smallskip
\bigskip

\bibliographystyle{siamplain}
\bibliography{DC-En}

\begin{thebibliography}{10}

\bibitem{bbo03}
{\sc I.~Babu$\breve{s}$ka, U.~Banerjee, and J.~E. Osborn}, {\em Survey of
  meshless and generalized finite element methods: a unified approach}, Acta
  Numerica,  (2003), pp.~1--125.

\bibitem{frangos2010surrogate}
{\sc L.~Biegler, G.~Biros, O.~Ghattas, M.~Heinkenschloss, D.~Keyes, B.~Mallick,
  Y.~Marzouk, L.~Tenorio, B.~V.~B. Waanders, and K.~Willcox}, {\em 7. Surrogate
  and Reduced-Order Modeling: A Comparison of Approaches for Large-Scale
  Statistical Inverse Problems}, John Wiley \& Sons, Ltd, 2010.

\bibitem{bilionis2013solution}
{\sc I.~Bilionis and N.~Zabaras}, {\em Solution of inverse problems with
  limited forward solver evaluations: a bayesian perspective}, Inverse
  Problems, 30 (2013), p.~015004.

\bibitem{brooks2011handbook}
{\sc S.~Brooks, A.~Gelman, G.~Jones, and X.-L. Meng}, {\em Handbook of markov
  chain monte carlo}, CRC press, 2011.

\bibitem{bui2013computational}
{\sc T.~Bui-Thanh, O.~Ghattas, J.~Martin, and G.~Stadler}, {\em A computational
  framework for infinite-dimensional bayesian inverse problems part i: The
  linearized case, with application to global seismic inversion}, SIAM Journal
  on Scientific Computing, 35 (2013), pp.~A2494--A2523.

\bibitem{burger2001level}
{\sc M.~Burger}, {\em A level set method for inverse problems}, Inverse
  problems, 17 (2001), p.~1327.

\bibitem{chan2004level}
{\sc T.~F. Chan and X.-C. Tai}, {\em Level set and total variation
  regularization for elliptic inverse problems with discontinuous
  coefficients}, Journal of Computational Physics, 193 (2004), pp.~40--66.

\bibitem{christen2005markov}
{\sc J.~A. Christen and C.~Fox}, {\em Markov chain monte carlo using an
  approximation}, Journal of Computational and Graphical statistics, 14 (2005),
  pp.~795--810.

\bibitem{chung2005electrical}
{\sc E.~T. Chung, T.~F. Chan, and X.-C. Tai}, {\em Electrical impedance
  tomography using level set representation and total variational
  regularization}, Journal of Computational Physics, 205 (2005), pp.~357--372.

\bibitem{cotter2013mcmc}
{\sc S.~L. Cotter, G.~O. Roberts, A.~M. Stuart, D.~White, et~al.}, {\em Mcmc
  methods for functions: modifying old algorithms to make them faster},
  Statistical Science, 28 (2013), pp.~424--446.

\bibitem{cui2016dimension}
{\sc T.~Cui, K.~J. Law, and Y.~M. Marzouk}, {\em Dimension-independent
  likelihood-informed mcmc}, Journal of Computational Physics, 304 (2016),
  pp.~109--137.

\bibitem{dunlop2017hierarchical}
{\sc M.~M. Dunlop, M.~A. Iglesias, and A.~M. Stuart}, {\em Hierarchical
  bayesian level set inversion}, Statistics and Computing, 27 (2017),
  pp.~1555--1584.

\bibitem{efendiev2005efficient}
{\sc Y.~Efendiev, A.~Datta-Gupta, V.~Ginting, X.~Ma, and B.~Mallick}, {\em An
  efficient two‐stage markov chain monte carlo method for dynamic data
  integration}, Water Resources Research, 41 (2005), p.~12423.

\bibitem{efendiev2013generalized}
{\sc Y.~Efendiev, J.~Galvis, and T.~Y. Hou}, {\em Generalized multiscale finite
  element methods (gmsfem)}, Journal of Computational Physics, 251 (2013),
  pp.~116--135.

\bibitem{efendiev2011multiscale}
{\sc Y.~Efendiev, J.~Galvis, and X.-H. Wu}, {\em Multiscale finite element
  methods for high-contrast problems using local spectral basis functions},
  Journal of Computational Physics, 230 (2011), pp.~937--955.

\bibitem{forrester2007multi}
{\sc A.~I. Forrester, A.~S{\'o}bester, and A.~J. Keane}, {\em Multi-fidelity
  optimization via surrogate modelling}, Proceedings of the royal society a:
  mathematical, physical and engineering sciences, 463 (2007), pp.~3251--3269.

\bibitem{gunzburger2017ensemble}
{\sc M.~Gunzburger, N.~Jiang, and M.~Schneier}, {\em An ensemble-proper
  orthogonal decomposition method for the nonstationary navier--stokes
  equations}, SIAM Journal on Numerical Analysis, 55 (2017), pp.~286--304.

\bibitem{hinze2005proper}
{\sc M.~Hinze and S.~Volkwein}, {\em Proper orthogonal decomposition surrogate
  models for nonlinear dynamical systems: Error estimates and suboptimal
  control}, in Dimension reduction of large-scale systems, Springer, 2005,
  pp.~261--306.

\bibitem{hou1997multiscale}
{\sc T.~Y. Hou and X.-H. Wu}, {\em A multiscale finite element method for
  elliptic problems in composite materials and porous media}, Journal of
  computational physics, 134 (1997), pp.~169--189.

\bibitem{jiang2018bayesian}
{\sc L.~Jiang and N.~Ou}, {\em Bayesian inference using intermediate
  distribution based on coarse multiscale model for time fractional diffusion
  equations}, Multiscale Modeling \& Simulation, 16 (2018), pp.~327--355.

\bibitem{jiang2014algorithm}
{\sc N.~Jiang and W.~Layton}, {\em An algorithm for fast calculation of flow
  ensembles}, Int. J. Uncertain. Quantif, 4 (2014), pp.~273--301.

\bibitem{lee2013bayesian}
{\sc J.~Lee and P.~Kitanidis}, {\em Bayesian inversion with total variation
  prior for discrete geologic structure identification}, Water Resources
  Research, 49 (2013), pp.~7658--7669.

\bibitem{li2017novel}
{\sc Q.~Li and L.~Jiang}, {\em A novel variable-separation method based on
  sparse and low rank representation for stochastic partial differential
  equations}, SIAM Journal on Scientific Computing, 39 (2017),
  pp.~A2879--A2910.

\bibitem{liu2008monte}
{\sc J.~S. Liu}, {\em Monte Carlo strategies in scientific computing}, Springer
  Science \& Business Media, 2008.

\bibitem{martin2012stochastic}
{\sc J.~Martin, L.~C. Wilcox, C.~Burstedde, and O.~Ghattas}, {\em A stochastic
  newton mcmc method for large-scale statistical inverse problems with
  application to seismic inversion}, SIAM Journal on Scientific Computing, 34
  (2012), pp.~A1460--A1487.

\bibitem{marzouk2009stochastic}
{\sc Y.~Marzouk and D.~Xiu}, {\em A stochastic collocation approach to bayesian
  inference in inverse problems}, Communications in Computational Physics, 6
  (2009), pp.~826--847.

\bibitem{marzouk2009dimensionality}
{\sc Y.~M. Marzouk and H.~N. Najm}, {\em Dimensionality reduction and
  polynomial chaos acceleration of bayesian inference in inverse problems},
  Journal of Computational Physics, 228 (2009), pp.~1862--1902.

\bibitem{noor1980reduced}
{\sc A.~K. Noor and J.~M. Peters}, {\em Reduced basis technique for nonlinear
  analysis of structures}, Aiaa journal, 18 (1980), pp.~455--462.

\bibitem{nouy2010proper}
{\sc A.~Nouy}, {\em Proper generalized decompositions and separated
  representations for the numerical solution of high dimensional stochastic
  problems}, Archives of Computational Methods in Engineering, 17 (2010),
  pp.~403--434.

\bibitem{osher1988fronts}
{\sc S.~Osher and J.~A. Sethian}, {\em Fronts propagating with
  curvature-dependent speed: algorithms based on hamilton-jacobi formulations},
  Journal of computational physics, 79 (1988), pp.~12--49.

\bibitem{ou2019new}
{\sc N.~Ou, L.~Jiang, and G.~Lin}, {\em A new bi-fidelity model reduction
  method for bayesian inverse problems}, International Journal for Numerical
  Methods in Engineering,  (2019).

\bibitem{peherstorfer2018survey}
{\sc B.~Peherstorfer, K.~Willcox, and M.~Gunzburger}, {\em Survey of
  multifidelity methods in uncertainty propagation, inference, and
  optimization}, SIAM Review, 60 (2018), pp.~550--591.

\bibitem{ras2006gaussian}
{\sc C.~E. Rasmussen and C.~K. Williams}, {\em Gaussian process for machine
  learning}, MIT press, 2006.

\bibitem{Robinson2008Surrogate}
{\sc T.~D. Robinson, M.~S. Eldred, K.~E. Willcox, and R.~Haimes}, {\em
  Surrogate-based optimization using multifidelity models with variable
  parameterization and corrected space mapping}, Aiaa Journal, 46 (2008),
  pp.~2814--2822.

\bibitem{rudin1992nonlinear}
{\sc L.~I. Rudin, S.~Osher, and E.~Fatemi}, {\em Nonlinear total variation
  based noise removal algorithms}, Physica D: nonlinear phenomena, 60 (1992),
  pp.~259--268.

\bibitem{stuart2010inverse}
{\sc A.~M. Stuart}, {\em Inverse problems: a bayesian perspective}, Acta
  Numerica, 19 (2010), pp.~451--559.

\bibitem{vogel2002computational}
{\sc C.~R. Vogel}, {\em Computational methods for inverse problems}, vol.~23,
  Siam, 2002.

\bibitem{xiu2002wiener}
{\sc D.~Xiu and G.~E. Karniadakis}, {\em The wiener--askey polynomial chaos for
  stochastic differential equations}, SIAM journal on scientific computing, 24
  (2002), pp.~619--644.

\bibitem{yan2015stochastic}
{\sc L.~Yan and L.~Guo}, {\em Stochastic collocation algorithms using
  $l1$-minimization for bayesian solution of inverse problems}, SIAM Journal on
  Scientific Computing, 37 (2015), pp.~A1410--A1435.

\bibitem{yan2017convergence}
{\sc L.~Yan and Y.-X. Zhang}, {\em Convergence analysis of surrogate-based
  methods for bayesian inverse problems}, Inverse Problems, 33 (2017),
  p.~125001.

\bibitem{yan2019adaptive}
{\sc L.~Yan and T.~Zhou}, {\em Adaptive multi-fidelity polynomial chaos
  approach to bayesian inference in inverse problems}, Journal of Computational
  Physics, 381 (2019), pp.~110--128.

\bibitem{yao2016tv}
{\sc Z.~Yao, Z.~Hu, and J.~Li}, {\em A tv-gaussian prior for
  infinite-dimensional bayesian inverse problems and its numerical
  implementations}, Inverse Problems, 32 (2016), p.~075006.

\bibitem{zhou2015weighted}
{\sc T.~Zhou, A.~Narayan, and D.~Xiu}, {\em Weighted discrete least-squares
  polynomial approximation using randomized quadratures}, Journal of
  Computational Physics, 298 (2015), pp.~787--800.

\bibitem{Zhu2014Computational}
{\sc X.~Zhu, A.~Narayan, and D.~Xiu}, {\em Computational aspects of stochastic
  collocation with multifidelity models}, Society for Industrial \& Applied
  Mathematics, 2 (2014), p.~444–463.

\end{thebibliography}

\end{document}